\documentclass[letterpaper,12pt,leqno,oneside]{article}
\usepackage{color}
\usepackage{bm}
\usepackage{exscale}
\usepackage{amsmath}
\usepackage{amsfonts}
\usepackage{wasysym}
\usepackage{stmaryrd}
\usepackage{amscd}
\usepackage{graphicx}
\usepackage{amsxtra}
\usepackage{amssymb}
\usepackage{theorem}
\usepackage[final]{epsfig}
\usepackage{eqnarray}
\usepackage{hyperref}
\usepackage{enumitem}
\newcommand{\subscript}[2]{$#1 _ #2$}
\newtheorem{proposition}{Proposition}[subsection]
\newtheorem{definition}[proposition]{Definition}

\newtheorem{lemma}[proposition]{Lemma}

{\theorembodyfont{\rmfamily}\newtheorem{remark}[proposition]{Remark}}
\newtheorem{theorem}[proposition]{Theorem}
\newtheorem{corollary}[proposition]{Corollary}

{\theorembodyfont{\rmfamily}}

\newfont{\abc}{cmtt10 scaled 1200}

\def\R{\mathbb{R}}

\def\Z{\mathbb{Z}}

\def\U{\mathbb{U}}
\def\D{\mathbb{D}}
\def\P{\mathbb{P}}

\def\P{\mathbb{P}}

\def\E{\mathbb{E}}
\def\U{\mathbb{U}}

\def\TP{\mathbb{TP}}

\def\I{\mathbb{I}}
\def\ve{\varepsilon}

\def\ra{\rightarrow}
\def\cs{\symbol{35}}
\def\p{\partial}
\def\qed{\hfill $\Box$ \\}
\def\mm{\mbox}
\def\v{= \emptyset}

\def\D{\mathbf{ID}}
\def\M{\mathbb{A}}

\def\bp{\langle A \rangle}

\def\si{$\mathcal{S}$}

\def\bp{\langle A \rangle}

\begin{document}

\vspace*{-0.3cm}

\begin{center}\Large{\bf{Potential Theory on Minimal Hypersurfaces II: Hardy Structures and Schr\"odinger Operators}}\\
\bigskip
\large{\bf{Joachim Lohkamp}\\
\medskip}

\end{center}

\noindent Mathematisches Institut, Universit\"at M\"unster, Einsteinstrasse 62, Germany\\
 {\small{\emph{e-mail: j.lohkamp@uni-muenster.de}}}

{\small {\center \tableofcontents}

\vspace{0.3cm}
%
%
%
\setcounter{section}{1}
\renewcommand{\thesubsection}{\thesection}
\subsection{Introduction} \label{intro}
In this paper we further develop the potential theory on generally singular area minimizing and, more generally, on almost  minimizing hypersurfaces we initiated  in the first part of this work~\cite{L1}. We focus on elliptic operators \emph{naturally} associate to such hypersurfaces, for instance, from the minimality  constraint of the hypersurface and the geometry of the ambient space. This analysis can be applied to investigate the ambient space. A famous example hereof is the (non)existence of positive scalar curvature metrics on compact manifolds which is closely related to the potential theory of the conformal Laplacian on its minimal hypersurfaces~\cite{KW}, ~\cite{SY}.\\

We know from ~\cite{L2} that the \emph{\si-uniformity} of $H \setminus \Sigma$, for an (almost) minimizing hypersurface $H$ with singular set $\Sigma$ suffices to  establish a transparent potential theory on $H\setminus\Sigma$.   This comes along with boundary Harnack inequalities $H \setminus \Sigma$ relative the boundary $\Sigma$ and a Martin theory saying that $\Sigma$ is just the Martin boundary for a large class of operators, the so-called \emph{\si-adapted operators}.\\

In the present paper we take the applications of \si-structures to the elliptic analysis on (almost) minimizers much further. We first show that (almost) minimizers actually carry an enhanced version of \si-structures, namely a so-called \emph{Hardy \si-structure}. These satisfy an additional geometric coupling relation between \si-structures and the curvature of almost minimizers. It follows that many classical elliptic operators, like the Jacobi field operator or the conformal Laplacian, are actually \si-adapted so that our techniques from \cite{L1} apply. Typically, these operators will be \emph{Schr\"odinger operators} $L=L(H)$ on $H \setminus \Sigma$, and the assignment $H \mapsto L(H)$ \emph{commutes} with blow-ups. That is, for infinite scalings $m \cdot H$ around any $p \in \Sigma \subset H$, $L(m \cdot H)$ converges to $L(C)$ on any tangent cone $C$ in $p$ as $m \ra \infty$. For these so-called \emph{natural} Schr\"odinger operators our theory works particularly nicely.\\

The question of understanding their ground states and Green's functions is of great importance. Towards the singularity set these characteristic solutions all have \emph{minimal growth} compared with other solutions. This property is the starting point for our further analysis. Namely, we prove the allegedly plausible but rather non-trivial fact that minimal growth of solutions $u>0$ on $H \setminus \Sigma$ of $L(H) \,f=0$ towards some singular point $p \in \Sigma$ \emph{persists under blow-ups}. \\

In particular we find that, for any tangent cone $C$ around $p$, solutions of $L(C) \,f=0$ which are induced from $u$ have minimal growth towards the singularity set $\Sigma_C$ of the cone. The Martin theory in  \cite{L1} shows that these solutions are unique up to constant multiples and they admit a separation of variables. From this observation we can start a scheme of inductive descent and continue with blow-ups in singular points of $C \setminus \{0\}$ to get simpler cones of the form $C^* \times \R$ until we reach $C^\circ \times \R^k$, where  $C^\circ$ is only singular at its tip. The resulting tree of blow-ups together with the induced solutions, analyzed using the individual potential theory on these blow-ups,  allows a detailed analysis of solutions on $H \setminus \Sigma$ near $\Sigma$.\\

As in ~\cite{L2} this analysis does not use the structure of $\Sigma$ itself. Instead it is naturally and tightly connected to \si-structures and the properties of the associated hyperbolic unfoldings.
%
\subsubsection{Basic Concepts and Notations}\label{notation}%
We recall some basic notations from the first part of this work, cf.\ \cite[Chapter 1.1]{L1} and \cite[Appendix A]{L2} for details.\\

In this paper $H^n$ denotes a connected integer multiplicity rectifiable current of dimension $n\ge 2$ inside some complete, smooth Riemannian manifold $(M^{n+1},g_M)$. By $\Sigma_H$ (or simply $\Sigma$ if there is no risk of confusion) we denote the set of singular points of $H$. It is known to have Hausdorff dimension $\le  \dim  H-7$. We also write $\sigma_C$, in place of $\Sigma_C$, for the singular set of a minimal cone $C$, in particular if we want to emphasize that $C$ is viewed as a tangent cone.\\

The induced Riemannian metric on $H$ will be denoted by $g_H$. For $\lambda>0$ we let $\lambda\cdot M$ or $\lambda\cdot H$ denote the conformally rescaled Riemannian manifolds $(M,\lambda^2 \cdot  g)$ or $(H,\lambda^2 \cdot g_H)$. We refer to the induced distance function $d=d_H$ on $H$ as \emph{intrinsic distance}. Further, $A=A_H$ is the second fundamental form of $H \subset M$ and $|A_H|$ its norm.\\

We shall consider the following basic subclasses of such currents.
\begin{description}
  \item[${\cal{H}}^c_n$:] $H^n \subset M^{n+1}$ is  compact locally mass minimizing without boundary.
  \item[${\cal{H}}^{\R}_n$:] $H^n \subset\R^{n+1}$ is a complete hypersurface in flat Euclidean space $(\R^{n+1},g_{eucl})$ with $0\in H$,
  that is an oriented minimal boundary of some open set in $\R^{n+1}$.
  \item[${\cal{H}}_n$:] ${\cal{H}}_n:= {\cal{H}}^c_n \cup {\cal{H}}^{\R}_n$ and ${\cal{H}} :=\bigcup_{n \ge 1} {\cal{H}}_n$. We briefly refer to $H \in {\cal{H}}$ as an \textbf{area minimizer}.
  \item[$\mathcal{C}_{n}$:] $\mathcal{C}_{n} \subset {\cal{H}}^{\R}_n$ is the space of area minimizing $n$-cones in $\R^{n+1}$ with tip in $0$.
  \item[$\mathcal{SC}_{n}$:] $\mathcal{SC}_n \subset\mathcal{C}_n$ is the subset of cones which are at least singular in $0$.
  \item[$\mathcal{K}_{n-1}$:] For any area minimizing cone $C \subset \R^{n+1}$ with tip $0$, we get the non-minimizing minimal hypersurface $S_C$ in the unit sphere
  \[
  S_C:= \p B_1(0) \cap C \subset S^n \subset  \R^{n+1}
  \]
  and we set ${\cal{K}}_{n-1}:= \{ S_C\,| \, C \in {\mathcal{C}_{n}}\}$. We write ${\cal{K}}= \bigcup_{n \ge 1} {\cal{K}}_{n-1}$ for the space of all such hypersurfaces $S_C$.
  \item[${\cal{G}}^c_n$:] $H^n \subset M^{n+1}$ is a compact embedded \emph{almost} minimizer with $\p H \v$. We set ${\cal{G}}^c :=\bigcup_{n \ge 1} {\cal{G}}^c_n.$
  \item [${\cal{G}}_n$:] ${\cal{G}}_n := {\cal{G}}^c_n \cup {\cal{H}}^{\R}_n$ and ${\cal{G}} :=\bigcup_{n \ge 1} {\cal{G}}_n$. These are the main classes considered in this paper.
\end{description}
We denote the one-point compactification of a hypersurface $H \in {\cal{H}}^{\R}_n$ by $\widehat{H}$. For the singular set $\Sigma_H$ of some $H \in {\cal{H}}^{\R}_n$ we \emph{always} add $\infty_H$   to $\Sigma$ as well, even when $\Sigma$ is already compact, to define $\widehat\Sigma_H:=\Sigma_H \cup \infty_H$. On the other hand, for $H \in{\cal{G}}^{c}_n$ we set $\widehat H=H$ and $\widehat\Sigma=\Sigma$.\\

We call an assignment $\bp$ which associates with any $H \in {\cal{G}}$ a function $\bp_H:H \setminus \Sigma_H\to\R$ an \textbf{\si-transform} provided it satisfies the subsequent list of axioms.
\emph{\begin{description}
  \item[(S1)] \emph{\textbf{Trivial Gauge}} \,
  If $H \subset M$ is totally geodesic, then $\bp_H \equiv 0$.
  \item[(S2)] \emph{\textbf{\si-Properties}} \,
  If $H$ is not totally geodesic, then the level sets $\M_c:= \bp_H^{-1}(c)$, for $c>0$, we call the $|A|$-\textbf{skins}, surround the level sets of $|A|$:
  \[
  \bp_H>0, \bp_H \ge |A_H| \mm{ and } \bp_H(x) \ra \infty, \mm{ for } x \ra p \in \Sigma_H.
  \]
Like $|A_H|$, $\bp_H$ anticommutes with scalings, i.e., $\bp_{\lambda \cdot H} \equiv \lambda^{-1} \cdot  \bp_{H}$ for any $\lambda >0$.
  \item[(S3)] \emph{\textbf{Lipschitz regularity}} \,
  If $H$ is not totally geodesic, and thus $\bp_H>0$, we define the \textbf{\si-distance}  $\delta_{\bp_H}:=1/\bp_H$.
  This function is
  $L_{\bp}$-Lipschitz regular for some constant $L_{\bp}=L(\bp,n)>0$, i.e.,
  \[
  |\delta_{\bp_H}(p)- \delta_{\bp_H}(q)|   \le L_{\bp} \cdot d_H(p,q) \mm{ for any } p,q \in  H \setminus \Sigma \mm{ and any } H \in {\cal{G}}_n.
  \]
If $H$ is totally geodesic, and thus $\bp_H \equiv 0$, we set $\delta_{\bp_H}\equiv\infty$ and  $|\delta_{\bp_H}(p)- \delta_{\bp_H}(q)|\equiv0$.
  \item[(S4)]  \emph{\textbf{Naturality}} \,
  If $H_i \in {\cal{H}}_n$, $i \ge 1$, is a sequence converging* to the limit space $H_\infty \in {\cal{H}}_n$, then $\bp_{H_i}\overset{C^\alpha}  \longrightarrow {\bp_{H_\infty}}$ for any $\alpha \in (0,1)$. For general $H \in {\cal{G}}_n$, this holds for blow-ups:  $\bp_{\tau_i \cdot H} \overset{C^\alpha}  \longrightarrow {\bp_{H_\infty}}$, for any sequence $\tau_i \ra \infty$ so that  $\tau_i \cdot H \ra H_\infty \in {\cal{H}}^\R_n$.\\
\end{description}}

\noindent\textbf{Remark} \,
1.\ A construction of a concrete \si-transform by merging the metric $g_H$ on $H \setminus \Sigma_H$ and the second fundamental form $A=A_H$ into one scalar function $\bp_H$ on $H \setminus\Sigma_H$ was given in \cite{L2}.

\noindent2.\ The \si-distance $\delta_{\bp}$ is merely Lipschitz regular, but admits a Whitney type $C^\infty$-smoothing $\delta_{\bp^*}$ which satisfies (S1)-(S3) and which is quasi-natural in the sense that $c_1 \cdot \delta_{\bp}(x) \le \delta_{\bp^*}(x)  \le c_2 \cdot \delta_{\bp}(x)$ for some constants  $c_1,\,c_2>0$, cf.\ \cite[Proposition B.3]{L2} for details.\qed

%
\subsubsection{Statement of Results}\label{1}
In \cite{L1} we introduced \si-adapted operators and studied their basic potential theory on almost minimizers $H \in {\cal{G}}$. This class of operators is rich but a priori unrelated to the geometry of $(H \setminus \Sigma,g_H)$. Further, it depends on the chosen \si-transform $\bp$. In the first part of this paper we resolve this (apparent) issue. We show that there are \si-transforms with a tighter coupling between $|A|$ and $\bp$ beyond the axioms (S1)-(S4). Then we use this to prove that many classical operators are \si-adapted. Typically, these operators are \emph{symmetric} and we are mostly interested in eigenvalue problems, so we shall focus on this situation. We recall the following definitions and results from \cite[Definition 1, Theorem 7 and 8]{L1}. For this we use special charts for $H \setminus \Sigma$, namely \emph{\si-adapted charts}. These are bi-Lipschitz charts $\psi_p:B_R(p)\to\R$ centered in $p\in H\setminus\Sigma$, for some Lipschitz constant independent of $p$, and where the radius $R$ of the ball is, up to some common constant, just $1/\bp_H(p)$, cf.\ \cite[Chapter 2.3]{L1} and \cite[Proposition B.1]{L2}.\\

\noindent\textbf{Definition} \,
\emph{Let $H \in \cal{G}$. We call a symmetric second order elliptic operator $L$ on $H
\setminus \Sigma$  \textbf{shifted \si-adapted} supposed the following two conditions hold:}

\medskip

\textbf{$\mathbf{\bp}$-Adaptedness}  \,
\emph{$L$ satisfies \si-weighted uniformity conditions with respect to the charts $\psi_p$. Namely, we can write
\[
-L(u) = \sum_{i,j} a_{ij} \cdot \frac{\p^2 u}{\p x_i \p x_j} + \sum_i b_i \cdot \frac{\p u}{\p x_i} + c \cdot u,
\]
\emph{for some locally} $\beta$-H\"{o}lder continuous coefficients $a_{ij}$, $\beta \in (0,1]$, measurable functions $b_i$ and $c$, and there exists a $k_L=k \ge 1$ such that for any $p \in H\setminus \Sigma$ and $\xi\in\R^n$:}
\begin{itemize}
  \item $k^{-1} \cdot\sum_i \xi_i^2 \le \sum_{i,j} a_{ij}(p)\cdot \xi_i \xi_j \le k \cdot \sum_i \xi_i^2$,
  \item $\delta^{\beta}_{\bp}(p)\cdot |a_{ij}|_{C^\beta(B_{\theta(p)}(p))} \le k$,
  \item $\delta_{\bp}(p)\cdot |b_i|_{L^\infty}\le k$ and $\delta^2_{\bp}(p) \cdot |c|_{L^\infty}\le k$.
\end{itemize}

\textbf{$\mathbf{\bp}$-Finiteness} \,
\emph{There exists a finite constant $\tau = \tau(L,\bp,H)>-\infty$ such that for any smooth function $f$ which is compactly supported in $H\setminus \Sigma$, we have}
\begin{equation}\label{hadi0}
\int_H  f  \cdot  L f  \,  dV \, \ge \, \tau \cdot \int_H \bp^2\cdot f^2 dV.
\end{equation}

\noindent\textbf{Definition} \,
\emph{The largest $\tau \in (-\infty,+\infty)$  such that \eqref{hadi0} holds is the \textbf{principal eigenvalue} $\lambda^{\bp}_{L,H}$ of $\delta_{\bp}^2 \cdot L$.  The operator $L$ is called \textbf{\si-adapted} if $\lambda^{\bp}_{L,H}>0$.}\\

For any shifted \si-adapted $L$ we set
\[
L_\lambda:= L - \lambda \cdot \bp^2 \cdot Id, \mm{ for }\lambda \in\R.
\]
Again, $L_\lambda$ is shifted \si-adapted. Moreover, it is \si-adapted if and only if $\lambda < \lambda^{\bp}_{L,H}$.\\

Coming back to the inequality $\bp_H \ge |A_H|$ we note that there cannot be a \emph{pointwise} inverse inequality. For instance, singular cones may contain subcones where $|A|\equiv 0$ whereas $\bp >0$. However, we can prove an inverse \emph{integral} inequality for special \si-transforms. To exclude the previous counterexample we need to include a gradient into the integrals. We obtain thus a generalization of both the Poincar\'{e} inequality and the sharper Hardy inequality for the Laplacian $-\Delta_{Eucl}$ on flat Euclidean domains $D \subset \R^n$, to operators which couple to $|A|$ on curved manifolds $H\setminus \Sigma$ with boundary $\Sigma$, $H \in {\cal{G}}$, cf.Ch.\ref{hardy}.\\

\noindent\textbf{Theorem 1} \textbf{(Hardy \si-Structures, see Theorem~\ref{hhh})} \,
\emph{There are \si-transforms $\bp$ such that for all $H \in {\cal{G}}$ and $f \in C_0^\infty(H \setminus \Sigma)$, the space of smooth functions compactly supported in $H \setminus \Sigma$, the following \textbf{Hardy relations} hold:}
\begin{description}
  \item[($H_c$)] \emph{For any compact $H \in {\cal{G}}^c$ and any smooth $(2,0)$-tensor $B$ on the ambient space $M$ of $H$ with $B|_H \not\equiv -A_H$ there exists a constant $k_{H;B} > 0$ such that}
  \[
  \int_H|\nabla f|^2 + |A+B|_H|^2 \cdot f^2 dV \ge k_{H;B}  \cdot \int_H \bp^2 \cdot f^2 dV \ge \frac{k_{H;B}}{L^2_{\bp}} \cdot \int_H \frac{f^2}{dist(x, \Sigma)^2}dV.
  \]
  \emph{For singular $H \in {\cal{G}}^c$, $|A|$ is unbounded whereas $|B|_H|$ remains bounded, hence, in this case the condition $B|_H \not\equiv -A$ is redundant.}
  \item[($H_\R$)] \emph{For $H \in{\cal{H}}^{\R}_n$ we only consider the case $B=0$. Then the Hardy constant depends solely on the dimension, that is, $k_{H;0}=k_n>0$, and we have}
  \[
  \int_H|\nabla f|^2 + |A|^2 \cdot f^2 dV \ge k_n  \cdot \int_H \bp^2 \cdot f^2 dV \ge   \frac{k_n}{L^2_{\bp}} \cdot \int_H \frac{f^2}{dist(x, \Sigma)^2} dV.
  \]
\end{description}
\emph{These relations also apply to the case $\Sigma \v$ via the convention $1/dist(x, \Sigma)=0$. An \si-transform satisfying both axioms $(H) = (H_c) + (H_\R)$ is called a \textbf{Hardy \si-transform}.}\\

\noindent\textbf{Remark} \,
1.\ The ambient field $B$ incorporates geometric or physical constraints on $M$. An example is the second fundamental form of $M$ in a still higher dimensional space. In this context, we consider hypersurfaces in ${\cal{H}}^{\R}_n$ primarily as limit spaces under blow-ups, that is, infinite scaling of a given  $H \in {\cal{G}}$ around its singularities. Since $|B|_{\lambda M}= \lambda ^{-1} \cdot |B|_{M}$, for  $\lambda >0$,  and $\lambda\ra \infty$ during the blow-up process,   this suggests to focus on the case $B=0$ on $H \in {\cal{H}}^{\R}_n$.

\noindent 2.\ ($H_c$) implies the Poincar\'{e} type inequality $\int_H|\nabla f|^2 + |A+B|_H|^2 \cdot f^2 dV \ge k^*_{H;B} \cdot \int_H  f^2 dV$ for some $k^*_{H;B} > 0$ since $\bp$ remains positively lower bounded when $H$ is compact. \qed\\

Henceforth, we only use \si-transforms on $\cal{G}$ which are Hardy. This allows us to show that many operators which occur in the Euler-Lagrange equations of natural variational integrals on $H$ are actually (shifted) \si-adapted. We have the following list of basic examples, where $scal_H$ and $scal_M$ denote the scalar curvature of $H$ and $M$ and $Ric_M$ the Ricci curvature of $M$.\\

\noindent\textbf{Theorem 2} \textbf{(Curvature Constraints, see Theorem~\ref{scal})} \,
\emph{Let $\bp$ be a Hardy \si-transform. Further, let $H^n=H \in \cal{G}$ with $H\subset M^{n+1}$ and such that $H \setminus \Sigma$ is non-compact and non-totally geodesic. The following operators are (shifted) \si-adapted on $H\setminus\Sigma$:
\begin{enumerate}
  \item The \textbf{conformal Laplacian} $L_H: = -\Delta_H +\frac{n-2}{4 (n-1)} \cdot scal_H$ is shifted \si-adapted. Furthermore, if $scal_M \ge 0$ and $H \in \cal{H}$, then $L_H$ is even \si-adapted.
  \item More generally, let $S$ be any smooth function on $M$. Then the \textbf{S-conformal Laplacian} $L_{H,S}: = -\Delta_H +\frac{n-2}{4 (n-1)} \cdot (scal_H - S|_H)$ is shifted \si-adapted. If $scal_M \ge S$ and $H \in \cal{H}$, then $L_{H,S}$ is even \si-adapted.
  \item The \textbf{Laplacian} $-\Delta_H$ is shifted \si-adapted. When $H$ is compact, the principal eigenvalue $\lambda^{\bp}_{-\Delta,H}$ vanishes and the ground state is given by a constant function. In particular, $H \setminus \Sigma$ has the Liouville property saying that all bounded harmonic functions are constant.
  \item For any smooth $(2,0)$-tensor $B$ on the ambient space $M$ with $B|_H \not\equiv -A$ if $H$ is compact and $B=0$ if $H \in {\cal{H}}^{\R}_n$, the \textbf{A+B Laplacian} $C_{H;A,B}:=-\Delta + |A+B|_H|^2$ is an \si-adapted operator.
  \item The \textbf{Jacobi field operator} $J_H=-\Delta_H - |A|^2-Ric_M(\nu,\nu)$ is shifted \si-adapted. Moreover, it has principal eigenvalue $\ge 0$ if $H \in \cal{H}$.
\end{enumerate}
}

These operators have two basic properties in common. First, they are \emph{naturally} associated with any $H \in {\cal{G}}$. This means that there is a unique expression for $L(H)$ on $H \setminus \Sigma_H$ such that the assignment $H \mapsto L(H)$ commutes with the convergence of sequences of almost minimizers, cf.\ Ch.~\ref{ndo} for details. Second, they are \emph{Schr\"odinger operators} with finite principal eigenvalues. We merge these properties into one concept (cf.\ also Chapters \ref{spli} - \ref{ivn}):\\

\noindent\textbf{Definition} \,
\emph{A natural and shifted \si-adapted operator $L$ is called a \textbf{natural Schr\"odinger operator} if for any given $H \in \cal{H}$, $L(H)$ has the form
\[
L(H)(u)=-\Delta_H \, u + V_H(x) \cdot \, u \mm{ on } H \setminus \Sigma_H
\]
for some H\"older continuous function $V_H(x)$.\\
}

We take the analysis of natural Schr\"odinger operators near singular points beyond the Martin theory on $H \setminus \Sigma_H$ by considering tangent cones and cone reduction arguments. These are common in the \emph{geometric} study of (almost) minimizers, for instance to prove bounds on the codimension of their singularity sets. For natural  Schr\"odinger operators we build a matching \emph{analytic} reduction scheme. In view of the key role played by the boundary Harnack inequalities ~\cite[Theorem 1 and 2]{L1} our goal is to understand how minimal growth towards singularities, transfers to the associated \emph{induced} solutions on the tangent cones obtained by blowing up.\\

The minimal growth concept we use in this context is that of solutions\emph{L}-vanishing in (parts of) $\Sigma_H$. We recall from ~\cite{L1} that a solution $u >0$  of $L \,f=0$ in $p \in \widehat{\Sigma}$  is \textbf{\emph{L}-vanishing} in $p$ when there exists a supersolution $w >0$ with $u/w(x) \ra 0$ for $x \ra p$ with $x \in H \setminus \Sigma$.\\

The following Theorem asserts an inheritance of minimal growth properties under blow-ups of the underlying spaces. The remarkable point is that we are comparing the fine asymptotic analysis of distinct spaces with completely different singularity sets. We are not aware of any comparable result in the literature.\\

\noindent\textbf{Theorem 3} \textbf{(Minimal Growth and Blow-Ups, see Theorem~\ref{miii})} \,
\emph{Let $H \in {\cal{G}}$ and $L$ be some natural Schr\"odinger operator on $H \setminus \Sigma$. Further, let $p \in \Sigma_H$ and $C$ be any tangent cone in $p$. If $L(H)$ is \si-adapted and $u >0$ a solution of $L(H) \, f =0$ which is $L$-vanishing in a neighborhood $V$ of $p$, then $L(C)$ is again \si-adapted. Further, any solution of $L(C) \, f =0$ induced by $u$ is $L(C)$-vanishing along the entire singular set $\sigma_C \subset C$.}\\

\noindent\textbf{Remark} \, If we do not fix the singular point while we scale $H$, that is, we consider a converging subsequence of pointed spaces $(s_i \cdot H, p_i)$ for $p_i \ra p$, $p_i \in \Sigma_{H}$, and $s_i \ra \infty$ of scaling factors, then the limit space $(H_\infty, p_\infty)$ can be a general Euclidean hypersurface in ${\cal{H}}^{\R}_n$. In this case any induced solution on $H_\infty \setminus \Sigma_{H_\infty}$ is $L$-vanishing along $\Sigma_{H_\infty}$. \qed\\

The \si-adaptedness of $L(C)$ implies in particular that Martin theory applies to $C \setminus \sigma_C$. By~\cite[Theorem 3]{L1} there is exactly one Martin boundary point $\Psi_+$ at $\infty$. We emphasize that it is not the symmetry of the cone but the \si-uniformity of $C$ which implies uniqueness of $\Psi_+$. Indeed, Ancona gave striking counterexamples of Euclidean cones over \emph{non-uniform} spherical domains with uncountable families of minimal Martin boundary points at infinity~\cite{An3}. In turn, the uniqueness of $\Psi_+$ \emph{combined} with the cone symmetry of $C$ and the separation of variables for natural Schr\"odinger operators over cones yields the following structure result.\\

\noindent\textbf{Theorem 4} \textbf{(Separation of Variables, see Theorem~\ref{fix1})} \,
\emph{Let $C \in \mathcal{SC}_{n}$ and $L$ be a natural Schr\"odinger operator. For the \si-adapted operator $L_\lambda=L - \lambda \cdot \bp^2 \cdot Id$, $\lambda < \lambda^{\bp}_{L,C}$, and the two distinguished points $\Psi_-$ at zero and $\Psi_+$ at infinity in the  Martin boundary of $L_\lambda(C)$, we have, in terms of polar coordinates $(\omega,r)$:
\[
\Psi_\pm(\omega,r) = \psi(\omega) \cdot r^{\alpha_\pm}, \, (\omega,r) \in S_C \setminus \sigma \times \R^{>0}, \mm{ with } \textstyle \alpha_\pm = - \frac{n-2}{2} \pm \sqrt{ \Big( \frac{n-2}{2} \Big)^2 + \mu}
\]
for some constant $\mu(C,L,\lambda)>-(n-2)^2/4$.}\\

These results describe the behavior of the $\Psi_\pm$ on individual cones. In general, however, we have infinitely many distinct tangent cones around a singular point $p \in \Sigma_H$. Our next result is a variant of Theorem 3 asserting that the assignment of $\Psi_\pm$ to the underlying cones is natural in the following sense.\\

\noindent\textbf{Theorem 5} \textbf{(Naturality of  $\Psi_\pm$, see Theorem~\ref{miv})}
\emph{ Let  $L$ be a natural and \si-adapted Schr\"odinger operator on cones $C \in   \mathcal{SC}_n$. Then for any flat norm converging sequence $C_i \ra C_\infty$, $i \ra \infty$, with suitably normalized associated solutions $\Psi_\pm(C_i)$ and $\Psi_\pm(C_\infty)$, we have
\[
\Psi_\pm(C_i) \circ \D \ra \Psi_\pm(C_\infty)\;\; C^{2,\alpha}\mm{-compactly on } C_\infty \setminus \sigma_{C_\infty} \mm{ as } i \ra \infty, ,
\]
where $\D$ is the asymptotic identification representing $C_i$ as a smooth section of the normal bundle over $C_\infty$, cf.\ Chapter \ref{tfr} and \cite[Chapter 1.3]{L2}.}\\

A more general  version of this result applies to converging sequences in ${\cal{H}}^{\R}_n$, cf. \ref{miiv}.\\

\noindent\textbf{Remark} \,
Theorems 3, 4 and 5 give us a recipe for the asymptotic analysis of natural \si-adapted Schr\"odinger operators near $\Sigma \subset H$. For instance, if $p\in \Sigma$ is an isolated point and its tangent cones are singular only in the tip, these results entail a sharp description of the solutions of minimal growth on $H \setminus \Sigma$ near $p$. On the other hand, we can treat more complicated singular sets with tangent cones singular also outside the tip as follows: We consider the tangent cone and an induced solution of minimal growth towards $\sigma_C \subset C$ as the new initial object. Now we blow-up in points of $\sigma_C$ distinct from the tip and iterate this process until we reach the elementary case of a product cone $\R^m \times C^{n-m} $, for some $C^{n-m} \subset \R^{n-m}$ singular only in $0$ (this happens at the latest after $\dim H -7$ times). Since by uniqueness, the induced minimal solutions are $\R^{n-m}$-translation symmetric we end up with the explicit description provided by Theorem 4 over cones singular only at the tip. Finally, Theorem 5 takes care of the non-uniqueness of tangent cones and yields uniform control for all these cones. We can then work backwards this tree of blow-ups to the initially given $H \in {\cal{H}}$ from such a terminal node $\R^m \times C^{n-m} $, for some $C^{n-m} \subset \R^{n-m}$ singular only in $0$. \qed

As an example (and with later applications in scalar curvature geometry and general relativity in mind) we derive more detailed results for the conformal Laplacians $L_H$ and $L_C$ on $H$ and its tangent cones $C$:\\

\noindent\textbf{Theorem 6} \textbf{(Conformal Laplacians, see Theorem~\ref{evee})} \,
\emph{There are constants $\Lambda_n > \lambda_n >0$ depending only on $n$ such that for $\lambda \in (0, \lambda_n]$ and any singular area minimizing cone $C$, $(L_C)_\lambda$ is \si-adapted. Furthermore, we have the following estimates for $\Psi_\pm(\omega,r) = \psi(\omega) \cdot r^{\alpha_\pm}$:}
\begin{itemize}
  \item $0 > \alpha_+ \ge - (1- \sqrt{\tfrac{3}{4}}) \cdot \tfrac{n-2}{2} > - \tfrac{n-2}{2} > - (1+ \sqrt{\tfrac{3}{4}}) \cdot \tfrac{n-2}{2} \ge \alpha_- > -(n-2)$
  \item $|\psi|_{L^1(S_C \setminus \Sigma_{S_C})} \le a_{n,\lambda} \cdot \inf_{\omega \in S_C \setminus \Sigma_{S_C}}\psi(\omega)$ \emph{for some constant $a_{n,\lambda} >0$ depending only on $n$ and $\lambda$.}
\end{itemize}

An interesting point in this result is the uniform separation of the exponents given by the lower bound $\alpha_+ - \alpha_- \ge  \sqrt{3/4} \cdot (n-2)$.
%
%
%
\setcounter{section}{2}
\renewcommand{\thesubsection}{\thesection}
\subsection{Hardy \si-Structures}
The goal of this chapter is to show that many classical operators are actually shifted \si-adapted with respect to a special subclass of \si-structures which we call \emph{Hardy \si-structures}. Similar to Martin theory for \si-adapted operators, this is rather a property of the underlying space than of the operators.

%
\subsubsection{Hardy Inequalities}\label{hardy}%
The Poincar\'{e} inequality is a frequently used tool in geometric analysis. For a Lipschitz regular and bounded Euclidean domain $D \subset \R^n$ it asserts that there is some constant $a_D > 0$ such that
\[
\int_D |\nabla \phi|^2 dV  \ge a_D \cdot \int_D   \phi^2 dV
\]
for any smooth function $\phi$ compactly supported in $D$. The Hardy inequality is a remarkable refinement of this result. Indeed, under the previous assumptions we can even find a constant $c_D > 0$ such that for any such $\phi$,
\[
\int_D |\nabla \phi|^2 dV  \ge c_D \cdot \int_D \frac{\phi^2}{dist(x,\p D)^2} dV,
\]
cf.\ the detailed expositions \cite{BEL} and \cite{GN} for some background information.\\

Now for $H \in {\cal{G}}$, the singular set $\Sigma \subset H$ plays the role of a \emph{boundary} of $H \setminus \Sigma$, and we shall prove similar Hardy inequalities as in the Euclidean case using the metric distance $dist(x,\Sigma)$. We refer to these as \emph{metric} Hardy inequalities. More importantly, we get for suitable \si-transforms a stronger \si-\emph{Hardy inequality} using the \si-distance $\delta_{\bp}(x)$. In fact, the metric Hardy inqualities are a simple byproduct of our \si-formalism, in the same way as uniformity of $H\setminus\Sigma$ follows from \si-uniformity~\cite{L2}.\\

We get two versions of Hardy inequalities according to whether $H$ is compact or lies in ${\cal{H}}^{\R}_n$. The latter case is the main case as it covers blow-up limits  of hypersurfaces in ${\cal{G}}$ around singular points.

\begin{theorem}[Hardy \si-Structures]\label{hhh}
There are \si-transforms $\bp$ such that for all $H \in {\cal{G}}$ and $f \in C_0^\infty(H \setminus \Sigma)$, the space of smooth functions compactly supported in $H \setminus \Sigma$, the following holds.
\begin{description}
  \item[($H_c$)] \emph{For any compact $H \in {\cal{G}}^c$ and any smooth $(2,0)$-tensor $B$ on the ambient space $M$ of $H$ with $B|_H \not\equiv -A_H$ there exists a constant $k_{H;B} > 0$ such that}
  \[
  \int_H|\nabla f|^2 + |A+B|_H|^2 \cdot f^2 dV \ge k_{H;B}  \cdot \int_H \bp^2 \cdot f^2 dV \ge \frac{k_{H;B}}{L^2_{\bp}} \cdot \int_H \frac{f^2}{dist(x, \Sigma)^2}dV.
  \]
  \emph{For singular $H \in {\cal{G}}^c$, $|A|$ is unbounded but $|B|_H|$ remains bounded, so that, in this case, the condition $B|_H \not\equiv -A$ is redundant.}
  \item[($H_\R$)] \emph{For $H \in{\cal{H}}^{\R}_n$ we only consider the case $B=0$. Then the Hardy constant depends solely on the dimension, that is, $k_{H;0}=k_n>0$, and we have}
  \[
  \int_H|\nabla f|^2 + |A|^2 \cdot f^2 dV \ge k_n  \cdot \int_H \bp^2 \cdot f^2 dV \ge   \frac{k_n}{L^2_{\bp}} \cdot \int_H \frac{f^2}{dist(x, \Sigma)^2} dV.
  \]
\end{description}
These relations also apply to the case $\Sigma \v$ via the convention $1/dist(x, \Sigma)=0$. An \si-transform satisfying both axioms $(H) = (H_c) + (H_\R)$ is called a \textbf{Hardy \si-transform}.
\end{theorem}

\begin{remark}\label{moti}
For compact $H \in {\cal{G}}$, we need both integrands $|\nabla f|^2$ and $|A+B|_H|^2$ on the left hand side of these inequalities. Indeed, since $H$ is compact and $codim\,\Sigma \ge 2$, the coarea formula (cf.\ for instance \cite[Theorem 2.1.5.3]{GMS}) shows that
\[
\inf_{f \in C_0^\infty(H \setminus \Sigma)} \int_H|\nabla f|^2 dV/\int_H \bp^2\cdot f^2 dV = 0.
\]
On the other hand, $|A|$ usually vanishes or at least converges to zero along suitable sequences of points approaching $\Sigma$, while $\bp$ converges to $+\infty$. (This corresponds to rays in the tangent cones along which the cones are totally geodesic). Then, for smooth functions supported around such points, we get
\[
\inf_{f \in C_0^\infty(H \setminus \Sigma)} \int_H  |A+B|_H|^2\cdot f^2 dV /\int_H \bp^2\cdot f^2 dV= 0.
\]
\qed
\end{remark}

Since the Hardy inequalities ($H_c$) and ($H_\R$) do not follow from the \si-axioms we have to revisit our construction of metric \si-transforms $\bp_{\alpha,H}$ which we briefly recall. If $H \in {\cal{G}}$ is totally geodesic, we set $\bp_{\alpha,H}\equiv 0$. Otherwise, for $c>0$ we set $\M_c(\alpha):=$ the boundary of the $\alpha/c$-distance tube $\U^\alpha_c $ of $|A|^{-1}[c,\infty)$ and define
\[
\bp_{\alpha,H}(x):= \sup \{c \,| \, x \in \overline{\U}^\alpha_c\},\;\; x \in H \setminus \Sigma.
\]
A key input in proving that $\bp_{\alpha}$ also satisfies $(H) = (H_c) + (H_\R)$ will be \cite[Lemma A.7]{L1}. This is, roughly speaking, a quantitative version of the fact that tangent cones in singular points are also singular.

\begin{proposition}\label{mets}
For any $\alpha >0$ the metric \si-transform $\bp_{\alpha}$ is a Hardy \si-transform.
\end{proposition}

\begin{remark}
Apart of the proof of Proposition~\ref{mets} our later applications shall only make use of the axioms (S1)-(S4) and (H), but not of any particular feature of $\bp_\alpha$. Once we have proven this result we therefore simply add (H) to our set of \si-axioms and henceforth assume that all our \si-transforms are Hardy \si-transforms. \qed
\end{remark}

For the proof of Proposition \ref{mets} we may assume that $H$ is not totally geodesic, whence $\bp_\alpha>0$. Otherweise $H$ is totally geodesic, so $\bp_\alpha(x)\equiv 0$, and the Hardy inequalities hold trivially. We also note that the largest constant $k=k_{\alpha,H}$ such that
\[
\int_H|\nabla f|^2 + |A+B|_H|^2 \cdot f^2 dV \ge k\cdot \int_H \bp_\alpha^2 \cdot f^2 dV
\]
holds for any $f \in C_0^\infty(H \setminus \Sigma)$ is nothing but the first eigenvalue $\lambda_{P_\alpha}$ of the weighted operator
\[
P_\alpha u := \bp_\alpha^{-2} \cdot(-\Delta u + |A+B|_H|^2 \cdot u).
\]
Since $\bp_\alpha$ is locally Lipschitz, elliptic theory shows that eigenfunctions of $P_\alpha$ are  $C^{2,\gamma}$-regular for any $\gamma \in (0,1)$, cf.\ \cite[Chapter 6.4]{GT}. Then $\lambda_{P_\alpha}$ can be written as a Rayleigh quotient
\begin{equation}\label{ray}
\lambda_{P_\alpha} = \inf_{f \in C^{\infty}_0(H \setminus \Sigma),\, f \not\equiv 0} \int_H |\nabla f|^2  + |A+B|_H|^2 \cdot f^2 dV \Big/\int_H \bp_\alpha^2\cdot f^2 dV.
\end{equation}

To estimate $\lambda_{P_\alpha}$ we localize the problem to Neumann eigenvalues on regular balls. Then we take covers  by such balls with controlled covering numbers and use them to derive a positive lower estimate for the eigenvalue $\lambda_{P_\alpha}$.

%
\subsubsection{Neumann Eigenvalues on Balls}\label{hoh}
The Neumann eigenvalues of $P_\alpha$ on balls are scaling invariant: For any ball $B_r(p) \subset H \setminus \Sigma$, $r > 0$ and scaling factor $\mu >0$ we have
\begin{equation}\label{eigsc}
\nu_\alpha(B_r(p)) = \nu_\alpha(\mu\cdot B_r(p)),
\end{equation}
for the integrands of the Rayleigh quotient \eqref{ray} $|\nabla f|^2 + |A+B|_H|^2$ (numerator) and $\bp^2$ (denominator) rescale by the same factor $\mu^{-2}> 0$.\\

Next we would like to establish a lower bound for  $\nu_\alpha(B_r(p))$. However, there is no uniform positive lower bound for $r \ra 0$ if $|A|(p) = 0$. Conversely, when the balls become too large, we can neither control their geometry nor understand the eigenvalues or the covering numbers.  On the other hand, when we approach $\Sigma$ we get better and better local approximations by \emph{Euclidean} hypersurfaces in ${\cal{H}}^{\R}_n$ (after rescaling to a unit size). They are singular and real analytic hypersurfaces in $\R^n$. This leads us to the idea to let $\bp_\alpha$ determine the radius of the balls when we are close to $\Sigma$. Towards this end we notice that $\bp_{\alpha}(x) \ge \alpha / dist(x,\Sigma)$ means that $B_{\alpha/\bp_\alpha(p)}(p) \cap \Sigma \v$ for any $p \in  H \setminus \Sigma$. For these balls we have the following estimates:

\begin{lemma}\label{rad}
For any $\alpha \in (0,1]$ and $\mu \in (0,1/2)$ there is a neighborhood $U_{\alpha, \mu}$ of $\Sigma$ as well as a constant $\zeta(U_{\alpha, \mu}) >0$ such that
\begin{equation}\label{rade}
\nu_\alpha(B_{\mu \cdot \alpha /\bp_\alpha(p)}(p))  \ge \zeta \mm{ for any } p \in U_{\alpha, \mu} \setminus \Sigma.
\end{equation}
\end{lemma}
\noindent\textbf{Proof} \,
To simplify notation we only consider the case $\alpha =1$. We will make explicit use of the definition of $\bp_1$. Let us assume that there is a sequence of points $p_i \in H \setminus
\Sigma$ with $p_i \ra p_\infty \in \Sigma$ and such that $\nu_1(B_{\mu/\bp(p_i)}(p_i)) \ra 0$. After scaling $B_{\mu/\bp_1(p_i)}$ by $\bp_1(p_i)$ to obtain $B_{\mu}(p_i) \subset
\bp_1(p_i) \cdot H$ we may assume that
\[
\sup \{|A|(x) \, | \, x \in B_1(p_i) \cap \bp_1(p_i) \cdot H\} = 1
\]
and
\[
\sup \{|B|(x) \, | \, x \in B_1(p_i) \cap \bp_1(p_i) \cdot H\} \le 1/i
\]
in view of the definition of $\bp$, the boundedness of $|B|$ and $|B|_{\lambda M}= \lambda ^{-1} \cdot |B|_{M}$, $\lambda >0$. Moreover, taking without loss of generality $p_\infty=0$, we get a subsequence of the pointed spaces $(\bp_1(p_i) \cdot H, p_i)$ which converge compactly in flat norm and thus in $C^k$-norm, $k \ge 0$, to $(H_\infty,0) \subset (\R^{n+1},0)$ by standard regularity theory. (We generally use $k = 5$ to control also second derivatives of curvatures.) We notice that
\[
\sup \{|A|(x) \, | \, x \in B_1(p_i) \cap \bp_1(p_i) \cdot H\} = 1
\]
and thus $\bp_1(0) = 1$. Here where we use the following inequality ~\cite[Lemma A.7]{L2}: For any $\mu \in (0,1]$, there is a constant $c(\mu, n)>0$ such that
\[
\sup \{|A|(x) \, | \, x \in B_{\mu}(0) \cap H_\infty\} \ge c.
\]
We claim that the first Neumann eigenvalue $\nu_1(B_{\mu}(0) \cap H_\infty)$ for $P_1$ on $B_{\mu}(0) \cap H_\infty$ is positive:
\[
\inf_{f \in C^\infty (B_{\mu}(0) \cap H_\infty), f \not\equiv 0} \int_{B_{\mu}(0) \cap H_\infty} |\nabla f|^2  + |A|^2 \cdot f^2 dV \Big/\int_{B_{\mu}(0) \cap
H_\infty} \bp_1^2\cdot f^2 dV > 0.
\]
Note that the tensor $B$ vanishes in the limit since as $\lambda \ra \infty$ for $\lambda ^{-1} \cdot |B|_{M} \ra 0$ in $C^k$-norm. Since there is a positive upper bound for $\bp_1$ over $B_{\mu}(0) \cap H_\infty$, it suffices to consider
the usual non-weighted Neumann eigenvalue
\[
\nu(B_{\mu}(0) \cap H_\infty) = \inf_{f \in C^{\infty}(B_{\mu}(0) \cap H_\infty),\, f \not\equiv 0}\int_{B_{\mu}(0) \cap H_\infty} |\nabla f|^2  + |A|^2 \cdot f^2 dV/\int_{B_{\mu}(0) \cap
H_\infty}  f^2 dV.
\]
Clearly, $\nu(B_{\mu}(0) \cap H_\infty) \ge 0$, and for $\nu(B_{\mu}(0) \cap H_\infty) = 0$, we would have a smooth positive function $u$ with $\Delta u = |A|^2 \cdot u$ and
vanishing normal derivative along $\p B_{\mu}(0) \cap H_\infty$. But then Stokes theorem would imply $\int_{B_{\mu}(0) \cap H_\infty} \Delta u =0$, whereas
$\int_{B_{\mu}(0) \cap H_\infty} |A|^2 \cdot u > 0$. Hence for sufficiently large $i$,
\[
\nu_1(B_{\mu/\bp_1(p_i)}(p_i)) \ge \nu_1(B_{\mu}(0) \cap H_\infty)/2 > 0
\]
which contradicts $\nu_1(B_{\mu/\bp_1(p_i)}(p)) \ra 0$. \qed

%
\subsubsection{Controlled Covers}\label{coc}
To derive Proposition~\ref{mets} we combine these estimates for Neumann eigenvalues with the following \si-adapted covers for $H \in {\cal{H}}$  introduced in~\cite[Proposition B.1]{L2} to construct the \si-Whitney smoothings (which equally apply to $H \in {\cal{G}}$ since the whole argument is based on blow-up arguments to limits in  ${\cal{H}}^\R$). For a consistent statement we also include the trivial case of totally geodesic $H \in {\cal{G}}$.

\begin{proposition}[\si-Adapted Covers]\label{besi}
For any $H \in {\cal{G}}$ and size parameter $\xi \in (0,\xi_0)$ for some  $\xi_0(n,L_{\bp_{\alpha}}) \in (0,1/(10^3 \cdot L_{\bp_{\alpha}}))$ we get a locally finite cover ${\cal{A}} = \{\overline{B_{\Theta(p)}}\,|\, p \in A \}$ of $H \setminus \Sigma$ by closed balls of radius $\Theta(p):=   \xi / \bp_{\alpha}(p) = \xi \cdot \delta_{\bp_{\alpha}}(p)$ and some discrete set
$A \subset H \setminus \Sigma$, such that for a suitably small neighborhood $Q$ of $\Sigma$ we have:
\begin{description}
  \item[(C1)] For  $p \in Q$ the exponential map $\exp_p|_{B_{100 \cdot \Theta(p)}(0)}$ is  bi-Lipschitz onto its image for some bi-Lipschitz constant $l(n) \ge 1$.
  \item[(C2)] $A^Q := A \cap Q$ splits into $c(n)$ disjoint families $A^Q(1),...,A^Q(c)$ with
  \begin{enumerate}
    \item $B_{10\Theta(p)}(p) \cap B_{10\Theta(q)}(q) \v$, for $p$ and $q$ in the same $A^Q(k)$,
    \item $q \notin \overline{B_{\Theta(p)}(p)}$, for any two $p, q \in A^Q$.
  \end{enumerate}
\end{description}
In particular, for any $z \in Q$ and $\rho \in (0,10)$ the covering number
\begin{equation}\label{conn}
\cs(A \cap Q,z,\rho):= \cs\{x \in A \cap Q \,|\, z \in B_{\rho \cdot \Theta(x)}(x)\} \le c(n) \end{equation}
by balls centered in $A \cap Q$ is uniformly bounded. Such a cover ${\cal{A}}$ will be called \textbf{\si-adapted}. These covers have the following properties.
\begin{enumerate}
  \item For $H \in  {\cal{H}}^\R$  we may choose $Q = H\setminus\Sigma$.
  \item For any $\ve  >0$,  there is a $\xi_{\ve} \in (0,\xi_0(n,L_{\bp}))$ such that for  every $p \in H \setminus \Sigma$ the exponential map $\exp_p|_{B_{100 \cdot \xi_{\ve}/ \bp_{\alpha}(p)}(0)}$ is bi-Lipschitz onto its image with bi-Lipschitz constant $1 + \ve$.
\end{enumerate}
\end{proposition}
\noindent{\bf Proof} of Proposition~\ref{mets} \,
We first consider ($H_c$). We start with the simplest case where $H$ is totally geodesic. Then $\bp_{\alpha} \equiv 0$ and the Hardy inequality becomes trivial. Thus we take $H$ not totally geodesic so that $\bp_{\alpha} > 0$ on $H$.\\

Assuming that $H$ is regular we have an upper bound for $\bp_{\alpha} > 0$ on $H$. It is therefore enough prove the positivity of the usual eigenvalue of $-\Delta + |A+B|_H|^2$. Now $|A+B|_H|^2 \ge 0$, and in some open set it is positive since $B|_H \not\equiv -A$. Thus, for any smooth positive function $u$  (including the first eigenfunction) we have  $\int_H\Delta u =0$, whereas $\int_H |A+B|_H|^2 \cdot u > 0$. Hence, the eigenvalue of $-\Delta + |A+B|_H|^2$ cannot be zero.\\

Now we turn to the main case where $H$ is singular. As above we notice that $\nu_{\alpha}(H \setminus \overline{W})> 0$ for sufficiently small and smoothly bounded neighborhoods $W$ of $\Sigma \subset H$. Indeed, $H \setminus \overline{W}$ contains a non-empty open ball where $|A+B|_H|>0$ if $W$ is sufficiently small, since there are $p_k \in H \setminus \Sigma$ where $|A|(p_k)\ge k$, for any $k \ge 1$, whereas $|B|_H|$ remains bounded. For fixed $W$ we have positive bounds for $\bp_{\alpha}$, and we can consider the standard Neumann eigenfunction. Again we can invoke Stokes' Theorem to infer that $\nu_{\alpha}(H \setminus \overline{W}) > 0$. We note in passing that once we found a neighborhood $W$ with $\nu_{\alpha}(H \setminus \overline{W}) > 0$, positivity continues to hold for all neighborhoods $W^*$ of $\Sigma$ with $W^* \subset W$. However,  this argument does not give a uniform positive lower bound while $W^*$ shrinks to $\Sigma$ since $\bp_{\alpha}$ diverges when we approach $\Sigma$.\\

Next we take an \si-adapted cover for $\xi=\mu \cdot \alpha$ and radii  $\Theta(p) =  \mu \cdot \alpha/ \bp_{\alpha}(p)$, and we choose $W \subset Q$, so that  ${\cal{A}}_Q = \{\overline{B_{\Theta(p)}}\,|\, p \in Q \}$
is covering of $W$ with covering number $c(n)$. To ease notation we will write $B(p)=B_{\mu \cdot \alpha/\bp_\alpha(p)}(p)$. Then for any $f \in C^{\infty}_0(H \setminus \Sigma), f \not\equiv 0$ we get the following estimate:
\begin{align*}
&\int_H |\nabla f|^2 +  |A+B|_H|^2 \cdot f^2 dV\\
\ge&\frac{1}{{c(n)}+1}\cdot \left(\int_{H \setminus \overline{W}} |\nabla f|^2  +  |A+B|_H|^2 \cdot f^2 \, dV +   \sum_{B(p) \in {\cal{A}}} \int_{B(p)} |\nabla f|^2  +  |A+B|_H|^2 \cdot f^2 \, dV\right)\\
\ge&\frac{1}{{c(n)}+1}\cdot \left( \nu_{\alpha}(H \setminus \overline{W}) \cdot \int_{H \setminus \overline{W}}\bp_{\alpha}^2 \cdot  |f|^2 \, dV + \sum_{B(p) \in
{\cal{A}}}\nu_{\alpha}(B(p)) \cdot \int_{B(p)} \bp_{\alpha}^2 \cdot  |f|^2 \, dV  \right)\\
\ge&\frac{1}{{c(n)}+1}\cdot \min \Big\{\nu_{\alpha}(H \setminus \overline{W}), \inf_{B(p) \in {\cal{A}}} \nu_{\alpha}(B(p))\Big\} \cdot \int_H \bp_{\alpha}^2 \cdot |f|^2 dV.
\end{align*}
On the other hand, the Neumann eigenvalues $\nu_{\alpha}(B(p))$ are uniformly bounded from below by $\nu_\alpha(B_{\mu \cdot \alpha /\bp_\alpha(p)}(p))  \ge \zeta$, cf.\ \eqref{rade} in Lemma~\ref{rad}, whence
\[
\lambda_{P_\alpha} \ge \frac{\min \{\nu_{\alpha}(H \setminus \overline{W}), \zeta\}}{{c(n)}+1}  > 0.
\]
Put differently the Hardy inequality holds for $k_{\alpha ,H} :=\lambda_{P_\alpha}$.\\

These arguments apply equally well in the case $H \in {\cal{H}}^\R$. In this case Proposition~\ref{besi} asserts that we get \si-adapted covers not only of small neighborhoods of $\Sigma$, but of the entire hypersurface. Hence the previous chain of inequalities now yields
\begin{align*}
&\int_H |\nabla f|^2  +  |A|^2 \cdot f^2 dV\\
\ge& \frac 1{c(n)} \cdot \sum_{B(p) \in {\cal{A}}} \int_{B(p)} |\nabla f|^2  +  |A|^2 \cdot f^2 \, dV\\
\ge&\frac 1{c(n)}  \cdot   \sum_{B(p) \in {\cal{A}}}\nu_{\alpha}(B(p)) \cdot \int_{B(p)} \bp_{\alpha}^2 \cdot  |f|^2 \, dV.
\end{align*}
Therefore, the Hardy inequality holds for any $H \in {\cal{H}}^\R$ and for $k_{\alpha ,n} :=\zeta/{c(n)}$.\qed

%
\subsubsection{Geometric Operators}\label{exp}
For the remainder of this paper we fix a Hardy \si-transform $\bp$ on $\cal{G}$. In this section we shall employ the Hardy axiom $(H) $ to verify the \si-adaptedness of some basic geometric operators.\\

The natural geometric operator
\[
C_{H;A,B}:=-\Delta + |A+B|_H|^2,
\]
which we call the \textbf{A+B Laplacian}, couples to the second fundamental form of $H\subset M$ and an additional $(2,0)$-tensor $B$ on $M$. Of course, this definition is suggested right from the definition of Hardy \si-transforms, and the operator will be useful to provide lower bounds for variational integrals involving other \si-adapted operators.

\begin{lemma}[A+B Laplacian]\label{base}
For any smooth $(2,0)$-tensor $B$ on $M$ with $B|_H \not\equiv -A$ if $H \in {\cal{G}}^c_n$ and with $B=0$ if $H \in {\cal{H}}^{\R}_n$, the A+B Laplacian $C_{H;A,B}$ is an \si-adapted operator.
\end{lemma}
\noindent\textbf{Proof} \,  Theorem~\ref{hhh} shows the $\bp$-finiteness of $C_{H;A,B}$ for a positive $\tau>0$. Thus we only need to verify the $\bp$-adaptedness. With respect to the charts $\psi_p$ we write
\[
-\Delta u + |A+B|_H|^2 \cdot u = \sum_{i,j} a_{ij} \cdot \frac{\p^2 u}{\p x_i \p x_j} + \sum_i b_i \cdot \frac{\p u}{\p x_i} + c \cdot u.
\]
We recall that in local coordinates the Laplacian  $\Delta \, u$ equals $\frac{1}{\sqrt{det \, g}} \cdot \sum_{i,j} \frac{\p}{\p x_i} \Big(\sqrt{det \, g} \cdot
g^{ij} \cdot \frac{\p \, u}{\p x_j}\Big)$ and that $|A+B|_H|^2$ is smooth. Hence the coefficients $a_{ij}$, $b_i$, $c$ are smooth. Since the charts $\psi_p$ are the geodesic coordinates around $p$ we have with respect to these charts $a_{ij}(p) = \delta_{ij}$, $b_i(p) = 0$ and $c(p)=|A+B|_H|^2(p)$. Moreover, $|B|_H|$  remains bounded. Thus, the unfolding correspondence \cite[Proposition 3.3]{L1} shows that $\delta_{\bp^*}^2 \cdot L$ is an adapted operator on $(H \setminus \Sigma, d_{\bp^*})$.   \qed

One source for \si-adapted operators are geometric and physical variational problems. Here, geometric properties of the ambient space like its curvature, or a given tensor $T$ coming from physical constraints, can translate into \si-adapted operators on the hypersurface.

\begin{theorem}[Curvature Constraints]\label{scal} \,
Let $H^n=H \in \cal{G}$ with $H\subset M^{n+1}$ and such that $H \setminus \Sigma$ is non-compact and non-totally geodesic. Let $Ric_M(\nu,\nu)$ denote the Ricci curvature of $M$ for a normal vector $\nu$ of $H$, and let $scal_H$ and $scal_M$ be the scalar curvature of $H$ and $M$.
\begin{enumerate}
  \item The \textbf{conformal Laplacian}
  \[
  L_H: = -\Delta_H +\frac{n-2}{4 (n-1)} \cdot scal_H
  \]
  is shifted \si-adapted. Furthermore, $L_H$ is \si-adapted if $scal_M \ge 0$ and $H \in \cal{H}$.
  \item More generally, let $S$ be any smooth function on $M$. Then the \textbf{$S$-conformal Laplacian}
  \[
  L_{H,S}: = -\Delta_H +\frac{n-2}{4 (n-1)} \cdot (scal_H - S|_H)
  \]
  is shifted \si-adapted. Furthermore, $L_{H,S}$ is \si-adapted if $scal_M \ge S$ and $H\in\cal{H}$.
  \item The Laplacian $-\Delta_H$ is shifted \si-adapted. If $H$ is compact, the principal eigenvalue $\lambda^{\bp}_{-\Delta,H}$ vanishes and the ground state is given by a constant function. In particular, $H \setminus \Sigma$ has the Liouville property saying that all bounded harmonic functions are constant.
  \item The \textbf{Jacobi field operator}
  \[
  J_H:=-\Delta_H - |A|^2-Ric_M(\nu,\nu)
  \]
  is shifted \si-adapted. Furthermore, it has principal eigenvalue $\ge 0$ if $H \in \cal{H}$.
\end{enumerate}
\end{theorem}
\noindent\textbf{Proof} \,
By the arguments of Lemma \ref{base}, the operators are adapted to $\bp$.\\

To show \si-adaptedness it remains to consider the validity of the Hardy inequality. Towards this end consider the Gau\ss-Codazzi equation
\begin{equation}\label{miiq}
|A_H|^2 + 2Ric_M(\nu,\nu)  =  scal_M - scal_H +(tr A_H)^2,
\end{equation}
where $tr A_H$ is the mean curvature of $H$. If $H$ is minimal, then $tr A_H = 0$. The terms $2Ric_M(\nu,\nu)$ and $scal_M$ remain bounded, whereas $|A_H|^2$ and $scal_H$ diverge when we approach $\Sigma$ on $H$. Also  $S|_H$ remains bounded. Therefore, $\bp \ge |A|$ shows that $\bp \ge k \cdot \big|scal_H - S|_H \big|$ for some $k=k(\bp, S) \ge 1$.\\

We now turn to the $S$-conformal Laplacian $L_{H,S}$. As an area minimizing hypersurface $H$ is also stable, that is, the second variation of the area functional is non-negative. Hence, if $f$ is a smooth function on $H$ with $supp \: f \subset H \setminus \Sigma$ and $\nu$ the outward
normal vector field over $H \setminus \Sigma$, then
\begin{align}
Area'(f) &=  \int_H tr A_H (z) \cdot f(z) \: dV = 0,\nonumber\\
Area''(f)&= \int_{H}|\nabla f|^2 - \left( |A|^2 + Ric_M(\nu,\nu) \right) \cdot f^2 \: dV \ge 0.\label{secv}
\end{align}
The Gau\ss-Codazzi equation \eqref{miiq} gives an equivalent formulation of $Area'' (f) \ge 0$, namely
\begin{align}
&\int_H  f  \cdot  L_H f  \,  dV\nonumber\\
= &\int_H | \nabla f |^2 + \frac{n-2}{4 (n-1)} \cdot scal_H  \cdot  f^2 \, d V\nonumber\\
\ge&\int_H \frac{n}{2 (n-1)} \cdot  |  \nabla f |^2 + \frac{n- 2}{4 (n-1)}\cdot \left( | A |^2 + scal_M \right) \cdot  f^2\,  dV.\label{kwsy}
\end{align}
Further, assuming $scal_M \ge S$ gives
\begin{align*}
&\int_H  f  \cdot  L_{H,S} f  \,  dV \\
=& \int_H | \nabla f |^2 + \frac{n-2}{4 (n-1)} \cdot (scal_H - S|_H)  \cdot  f^2 \, d V\\
\ge&\int_H \frac{n}{2 (n-1)} \cdot  |  \nabla f |^2 + \frac{n- 2}{4 (n-1)}\cdot \left( | A |^2 + scal_M - S|_H \right) \cdot  f^2\,  d V\\
\ge&\frac{n- 2}{4 (n-1)}\cdot \int_H f  \cdot  C_{H;A,0} f   \, d V\\
\ge& \tau(\bp,H)  \cdot \frac{n- 2}{4 (n-1)}\cdot \int_H \bp^2 \cdot f^2   \, d V
\end{align*}
for $\tau(\bp,H) >0$. \\

For the Laplacian the condition $\lambda^{\bp}_{-\Delta,H} >-\infty$ is obvious since the variational integral is just the Dirichlet integral $\int_H | \nabla f |^2   \, d V \ge 0$. If $H$ is compact, it is easy to see that the principal eigenvalue equals $0$. Here we apply again the coarea formula, cf.~\cite[Theorem 2.1.5.3]{GMS}, and use the fact that the codimension of $\Sigma$ is $\ge 2$. Every constant function $v$ solves $\Delta \, v=0$. Thus $v \equiv 1$ can be taken as the ground state for compact $H$. Moreover, $v:=u + \inf_{H \setminus \Sigma} u +1$ is a positive harmonic function for any bounded harmonic function $u$. Hence $v$, and therefore $u$, are constant functions.\\

Finally, we consider the Jacobi field operator $J_H$. From the minimality of $H$ and \eqref{secv},
\[
\int_{H} f \cdot J_H f\, dV  = \int_{H}|\nabla_H f|^2  -  \left( |A|^2 + Ric_M(\nu,\nu)  \right) \cdot f^2 \, dV \ge 0
\]
for any smooth function $f$ on $H$ with $supp \: f \subset H \setminus \Sigma$. \qed

%
%
%
\setcounter{section}{3}
\renewcommand{\thesubsection}{\thesection}
\subsection{Blow-Up Martin Theory}
The structure of the singularity set of almost minimizers does not directly transfer to their tangent cones. However, for \si-adapted operators which are naturally associated with minimizers we do have a non-trivial relationship between their Martin theories on the minimizer and their tangent cones. In particular, by our main result the minimal growth properties of solutions towards singularities are inherited by all tangent cones. This plausible result uses virtually any structural detail of our theory, from \si-uniformity of almost minimizers to boundary Harnack inequalities of \si-adapted operators.

%
\subsubsection{$\D$-Maps  and Tangent Cones} \label{tfr}
For the convenience of the reader we briefly recall some concepts, notations and results from~\cite[Chapter 1.3]{L2} and \cite[Chapter 3.1]{L1}. These are needed to compare analytic data on almost minimizers which are spatially close.\\

\noindent\textbf{$\D$-maps and Naturality~\cite[Ch.1.3]{L2}} \,
The general setup is this. Consider a sequence of pairs $H^n_i \subset M^{n+1}_i$, where the $M^{n+1}_i$ are complete manifolds and the $H^n_i$ are almost minimizers with basepoint $p_i \in H_i$. Furthermore, we assume that as $i \ra \infty$, the pointed manifolds $(M_i,p_i)$ are compactly $C^k$-converging, for some $k \ge 5$, to a pointed limit manifold $(M,p)$. Hence, for any given $R >0$ and sufficiently large $i \ge i_R$ we have diffeomorphisms
\[
\Psi_i : B_R(p_i) \ra B_R(p) \mm{ such that } |(\Psi_i)_*(g_{M_i}) - g_{M}|_{C^k} \ra 0 \mm{ on } B_R(p).
\]
We assume that $\p H_i \cap B_R(p_i) \v$.
\begin{itemize}
  \item Basic compactness results show that the $H_i$ subconverge in flat norm to a limit area minimizer $H \subset M$ containing $p$, that is, as $i \ra \infty$ the sequence $\Psi_i (B_R(p_i) \cap H_i)$ subconverges to $B_R(p) \cap H$ in flat norm.
\end{itemize}
For instance, consider an initial hypersurface $H_0$ in $M_0$. A sequence $\tau_i \ra \infty$ gives rise to the rescaled sequence $H_i:=\tau_i \cdot H_0\subset M_i:=\tau_i\cdot M_0$ of blow-ups. Fix an $H_0$-singular point $p_0\in\Sigma\subset H_0$ and set $p_i:=p_0$. Then $M_i$ converges compactly to $\R^{n+1}$ while $H_i$ subconverges in flat norm to a (n area minimizing) tangent cone $H \subset \R^{n+1}$.
\begin{itemize}
  \item When $B_R(p) \cap H$ is \emph{smooth}, standard regularity shows the following. First, flat norm convergence of $\Psi_i(B_R(p_i) \cap H_i)$ to $B_r(p)\cap H$ implies that $B_R(p_i) \cap H_i$ is also smooth for sufficiently large $i$. Further, we obtain $C^k$-convergence in the following sense: The $\Psi_i (B_R(p_i) \cap H_i)$ can be identified with local $C^k$-sections $\Gamma_i :B_R(p) \cap H \ra \Psi_i (B_R(p_i) \cap H_i)\subset\nu$ of  the normal bundle $\nu$ of $B_R(p)\cap H$ in $M$, and $\Gamma_i$ $C^k$-converges to the zero section which we identify with $B_R(p) \cap H$.
\end{itemize}

\begin{definition}[$\D$-map]\label{idmap}
For sufficiently large $i$ we call the uniquely determined section
\[
\D := \Gamma_i : B_R(p) \cap H \ra  \Psi_i(B_R(p_i) \cap H_i)\subset\nu
\]
the \textbf{asymptotic identification} or $\D$\textbf{-map} for short.
\end{definition}

Put differently, the regularity theory of area minimizers implies that we can approximate smooth regions of $H$ by corresponding (and also smooth) regions of $H_i$  through the canonical local diffeomorphisms given by $\D$-maps, cf.\ also \cite[Chapter 1.3]{L2}.  Moreover, $\D$-maps are arbitrarily close to isometries in $C^k$-topology for sufficiently large $i$.  $\D$-maps are useful to visualize situations where we compare functions or operators on different underlying hypersurfaces in ${\cal{G}}$, for instance when passing to blow-up limits.

\begin{definition}[Natural assignements]\label{cr}
An assignment $F:H \mapsto F_H$, $H \in {\cal{G}}$, of functions $F_H:H \setminus\Sigma_H\to\R$ is called \textbf{natural} if $F_H$ commutes with the convergence of underlying spaces. More precisely, for any pointed sequence $H_i\in {\cal{G}}$, $p_i\in H_i\setminus\Sigma_{H_i}$ which is locally converging in flat norm to a pointed space $H$, $p\in H \setminus \Sigma_H$, there is a neighborhood $U(p) \subset H \setminus \Sigma_{H}$ such that
\[
|\D^*F_{H_i} - F_{H}|_{C^k(U(p))}=|F_{H_i} \circ \D - F_{H}|_{C^k(U(p))} \ra 0 \mm{ as } i \ra \infty
\]
for some $k=k(F) >0$.
\end{definition}

This naturality concept readily extends to more general assignments $H \mapsto F_H$ of tensors or operators defined over $H\setminus\Sigma$. In the latter case we refer to the operator $F(H)$ itself as \emph{natural}. It becomes an important problem to determine which properties of natural operators are stable under convergence, for instance under blowing up. As a first example we prove that (shifted) \si-adaptedness is stable.

\begin{lemma}[Inherited \si-Adaptedness]\label{inhnat}
Let $L$ be a natural operator and  $H \in {\cal{G}}$. Moreover, assume that $L(H)$ is (shifted) \si-adapted. Then $L(N)$ is (shifted) \si-adapted on any blow-up $N$ of $H$ around a singular point in $\Sigma \subset H$ and the principal eigenvalues satisfy $\lambda^{\bp}_{L,H} \le \lambda^{\bp}_{L,N}$.
\end{lemma}
\noindent\textbf{Proof} \,
For any smooth function $f$ with compact support $K \subset N \setminus \Sigma_N$ we choose a sufficiently large scaling factor $\gamma \gg 1$ such that the $\D$-map between $\gamma \cdot H$ and $N$ in a neighborhood of $K$ is very close to an isometry in $C^5$-topology. Then the naturality of $L$ and $\bp$ implies that
\[
\int_Hf \circ \D^{-1} \cdot L(\gamma\cdot H) (f\circ \D^{-1}) \, dV/ \int_H \bp_H^2  \cdot  f^2\circ \D^{-1} dV
\]
is arbitrarily close to
\[
\int_N  f \cdot L(N) (f) \,  dV/ \int_N \bp_N^2  \cdot  f^2  dV
\]
upon choosing $\gamma$ large enough. Thus the eigenvalue of $\delta_{\bp}^2 \cdot L(H)$ on $H$ is a lower estimate for the second integral. In particular, the Hardy inequality holds for $L(N)$ on $N$. The adaptedness property of $L(N)$ follows from the scaling invariance of the estimates and the adaptedness of $L(H)$. \qed\\

\noindent\textbf{Tangent Cone Freezing~\cite[Ch.4.1]{L1}} \,
Let $H \in {\cal{G}}$. Blowing up around $p\in \Sigma_H$ yields only converging \emph{sub}sequences. In particular, there is no canonical tangent cone which could serve as a local model. Nevertheless, we can approximate the local geometry near $z$ using tangent cones. For this, we take $\omega>0$ and consider the \textbf{\si-pencil pointing to} $p\in \Sigma$ defined by
\[
\P(p,\omega)=\P_H(p,\omega)= \{x \in H \setminus \Sigma_H \,|\, \delta_{\bp_H} (x) > \omega \cdot d_H(x,p)\}.
\]
In view of scaling arguments it is also useful to consider the \textbf{truncated \si-pencil}
\[
\TP(p,\omega,R,r)=\TP_H(p,\omega,R,r):=B_{R}(p) \setminus B_{r} (p) \cap  \P(p,\omega)\subset H.
\]
While we zoom into some given singular point $p \in H$, formally accomplished by scaling with increasingly large $\tau >0$, we observe that $\tau \cdot  \TP_H(p,\omega,R/\tau,r/\tau)$ is better and better $C^k$-approximated by the corresponding  truncated \si-pencil in some (usually changing) tangent cone. Intuitively, the twisting of $\P(p,\omega)$ slows down as $\tau\to\infty$ until it \emph{asymptotically freezes} and ressembles a cone-like geometry.

\begin{proposition}[Asymptotic \si-Freezing]\label{freez}
Let $H \in {\cal{G}}$ and $p \in \Sigma_H$. Further, pick an $\ve > 0$ and a pair $R > 1 > r >0$. Then for any given $\omega \in (0,1)$ and $k \in \Z^{\ge 0}$, we have a $\tau_{\ve , R , r, \omega,p,k} \ge1$ with the following property. \\

For any $\tau\ge \tau_{\ve , R , r, \omega,p,k}$ there is a tangent cone $C_p^\tau$ of $H$ in $p$ such that the rescaled truncated \si-pencil $\tau \cdot  \TP_H(p,\omega,R/\tau,r/\tau)= \tau \cdot  \left(B_{R/\tau}(p) \setminus B_{r/\tau} (p) \cap  \P_H(p,\omega)\right)\subset\tau\cdot H$ can be written as a smooth section $\Gamma_\tau$ with $|\Gamma_\tau|_{C^k} < \ve$ of the normal bundle over $\TP_{C^\tau_p}(0,\omega,R,r)\subset C_p^\tau$.
\end{proposition}

This can be concisely expressed in terms of $\D$-maps
\[
\D=\D^{\tau, R, r, \omega} : \TP_{C^\tau_p}(0,\omega,R,r) \ra \tau \cdot  \TP_H(p,\omega,R/\tau,r/\tau)
\]
satisfying
\[
|\D - id_{C^\tau_p}|_{C^5(\TP(0,\omega,R,r))} \le \ve
\]
for the norm on sections of the normal bundle, thinking of $id_{C^\tau_p}$ as the zero section of the normal bundle over ${C^\tau_p}$.

%
%
\subsubsection{Induced Solutions}\label{ndo}
The \si-freezing of Proposition \ref{freez} shows that the $\D$-maps are almost isometric maps between suitable domains with compact closure in $H\setminus\Sigma_{H}$ and arbitrarily large truncated pencils $\TP_{C^\tau_p}(0,\omega,R,r) \subset C_p^\tau$, where $R>0$ is arbitrarily large and $r$, $\omega>0$ is arbitrarily small.\\

Next we consider a natural operator $L$. Restricted to these truncated pencils, the coefficients of $L(C_p^\tau)$ are arbitrarily close in H\"older norm to those of the $\D$-pull-back of $L(H)$. Thus, for a positive solution $u$ of $L(H)\, f=0$, $u \circ \D$ solves an elliptic equation whose coefficients are arbitrarily close in H\"older norm to those of $L(H)$. Since weak solutions of $L\, f = 0$ are $C^{2,\alpha}$-regular, we get locally \emph{uniform} constants for elliptic regularity estimates and Harnack inequalities for positive solutions, cf.\ \cite{BJS}, \cite{E} or \cite{GT}.\\

Now we take a sequence $\tau_i \ra \infty$ as $i \ra \infty$. The embedding $C^{2,\beta} \subset C^{2,\alpha}$ is \emph{compact} on bounded domains for $\beta \in (0,\alpha)$. Hence, after normalizing the value of the locally defined almost-solutions $u \circ \D$ of $L(C_p^{\tau_i }) \, f = 0$ in some basepoint of $\TP_{C_p^{\tau_i }}(0,\omega,R,r) \subset C_p^{\tau_i }$, there is a subsequence compactly $C^{2,\beta}$-converging to a positive solution on $\TP_{C_p^{\tau_\infty}}(0,\omega,R,r) \subset C_p^{\tau_\infty}$ as $i \ra \infty$. This also uses Harnack inequalities away from the basepoint.\\

We observe that, while $i \ra \infty$, the $\D$-map becomes more and more isometric on ever-growing truncated pencils $\TP_{C_p^{\tau_i }}(0,\omega,R,r) \subset C_p^{\tau_i }$ as $R\to\infty$ and $r\to0$, eventually exhausting $C_p^{\tau_i}\setminus \sigma_{C_p^{\tau_i }}$ in the limit. Therefore, possibly after selecting further subsequences, this process induces a positive $C^{2,\beta}$-regular function $v$ solving $L(C_p^{\tau_i }\setminus \sigma_{C_p^{\tau_i }}) \, f = 0$  on the whole of $C_p^{\tau_i }\setminus \sigma_{C_p^{\tau_i }}$, that is, we obtain an \emph{entire} solution. This gives rise to the following functional version of Proposition \ref{freez}.

\begin{proposition}[Asymptotic Functional Freezing]\label{freezfu}
Let $H\in\cal{G}$ and $p \in \Sigma_H$. Further, let $\alpha\in(0,1)$, $\ve$, $\omega >0$ and $R \gg 1 \gg r >0$ with $R$ and $r$ sufficiently large respectively small. Then there exists $\tau^*(L,\ve,\omega, R , r, p)>0$ such that for any $\tau\ge \tau^*$ and any entire solution $u>0$ of $L(H) \, f = 0$ the following is true. There exists a tangent cone $C^\tau_p$ such that on $\TP_{C^\tau_p}(0,\omega,R,r)$ we have an entire solution $v>0$ of $L(C^\tau_p) \, f = 0$ with
\begin{equation}\label{ffefua}
|\D - id_{C}|_{C^5(\TP_{C^\tau_p}(0,\omega,R,r))} \le \ve\quad\mm{and}\quad|u \circ \D / v-1|_{C^{2,\alpha}(\TP_{C^\tau_p}(0,\omega,R,r) )} \le \ve.
\end{equation}
\end{proposition}
\noindent\textbf{Proof} \,
The first assertion is merely Proposition~\ref{freez}. Next, assume that for some $\ve >0$, $R>r>0$, $\omega>0$, there exists a sequence $\tau_i \ra \infty$ such that for any entire solution $v>0$ of  $L(C^{\tau_i}_p) \, f = 0$ on $C^{\tau_i}_p \setminus \sigma_{C^{\tau_i}_p}$ we would have
\[
|u \circ \D /v-1|_{C^{2,\beta}(\TP_{C^{\tau_i}_p}(0,\omega,R,r))}  \ge \ve.
\]
Compactness of the space of tangent cones ${\cal T}_p$ in $p$ gives a subsequence $\tau_{i_k}$ such that $C^{\tau_{i_k}}_p$ converges to some tangent cone $C_p$. Hence for any entire solution $v>0$ of  $L(C_p) \, f = 0$ on $C_p \setminus \sigma_{C_p}$ we would also have that
\begin{equation}\label{li}
|u \circ \D /v-1|_{C^{2,\beta}(\TP_{C_p}(0,\omega,R,r))}  \ge \ve.
\end{equation}
However, our discussion before the proposition shows that we can choose another subsequence $\tau_{i_{k_m}}$ of $\tau_{i_k}$ such that the normalized $u \circ \D$ induces a positive entire solution on $C_p$ contradicting \eqref{li}. \qed

\begin{corollary}[Asymptotic Functional Freezing on Cone Spaces]\label{freezfu2}
Let $\alpha\in(0,1)$, $\ve$, $\omega >0$ and $R \gg 1 \gg r >0$ with sufficiently large $R$ and small $r$. Then there exists a $\zeta(L,\ve,\omega , R , r)>0$ such that following holds. For any $C$ and $C' \in {\cal{C}}_n$ with $d_H(S_C,S_{C'}) \le \zeta$ we have
\[
|\D - id_{C}|_{C^5(\TP_{C}(0,\omega,R,r))} \le \ve,
\]
and any entire solution $u_{C'} >0$ of $L(C') \, f = 0$ induces an entire solution $u_{C,C'}>0$ of $L(C) \, f = 0$ with
\begin{equation}\label{ffefu}
|u_{C'}  \circ \D /u_{C,C'} -1|_{C^{2,\alpha}(\TP_{C}(0,\omega,R,r))}  \le \ve.
\end{equation}
\end{corollary}

Completely similarly we get such entire solutions on more general blow-up hypersurfaces in ${\cal{H}}^{\R}_n$ obtained from scaling an almost minimizer $H \in \cal{G}$ without fixing a basepoint.

\begin{definition}[Induced Solutions] \label{ind}
We call such a positive entire solution on the limit hypersurfaces in ${\cal{H}}^{\R}_n$ an \textbf{induced solution}.
\end{definition}

Since induced solutions are positive entire solutions we can employ the full Martin theory on spaces in ${\cal{H}}^{\R}_n$ and their \si-adapted operators to analyze them.

%
%
\subsubsection{Minimal Growth Stability} \label{rec}

Let $L$ be a natural \si-adapted Schr\"odinger operator $L$. We show that $L$-vanishing is a \emph{stable} property which persists under deformations of the underlying spaces. A far reaching consequence is that any positive solution of $L \, f=0$ which is $L$-vanishing around some singular $p \in \Sigma_H$ induces solutions on any tangent cone $C$ at $p$ which are $L$-vanishing towards the entire singular set $\sigma_C$. In other words, the asymptotic minimal growth of solutions and the convergence of the underlying spaces are interchangeable. This is a remarkable phenomenon as the singular set of the underlying space, and thus the meaning of $L$-vanishing, may change dramatically under deformations.

\begin{theorem}[Blow-Up Stability of $L$-Vanishing Properties]\label{miii}
Let $H \in {\cal{G}}$ and $L$ be some natural, \si-adapted Schr\"odinger operator on $H \setminus \Sigma$. Further, let $u >0$ be a solution of $L \, f =0$ which is $L$-vanishing in a neighborhood $V$ of some point $p \in \Sigma$. Then we have the following inheritance results:
\begin{itemize}
  \item If $C$ is a tangent cone of $H$ in $p$, then any solution induced on $C$ is $L$-vanishing along $\sigma_C$.
  \item More generally, for a sequence $s_i \ra \infty$ of scaling factors and a sequence of points $p_i \ra p$ in $\Sigma_{H}$ such that $(s_i \cdot H, p_i)$ subconverges to a limit space $(H_\infty, p_\infty)$ with $H_\infty \in {\cal{H}}^{\R}_n$, any induced solution on $H_\infty \setminus \Sigma_{H_\infty}$ $L$-vanishes along $\Sigma_\infty$.
\end{itemize}
\end{theorem}

The Martin theory for these limit spaces says that there is, up to multiples, precisely one positive solution $L$-vanishing along the singular set $\Sigma_{H_\infty}$. It is  the unique and minimal Martin boundary point at infinity $\infty _{H_\infty}$ cf.\cite[Theorem 3]{L1}. We write this solution as  $\Psi_+= \Psi_+(H_\infty,L)$. Similarly, for any cone $C \in \mathcal{SC}_n$, there is also a unique (minimal) Martin boundary point at the origin $0_C$, the positive solution $\Psi_-= \Psi_-(C,L)$ which is $L$-vanishing along $\sigma_C \cup \{\infty_C\} \setminus \{0_C\}$.\\

A basic problem with tangent cones in singular points is that they are generally non-unique. The following second stability theorem partial compensates this issue. We show that the induced solutions change \emph{continuously} when we move from one to another tangent cone. Therefore we can extend compactness results for area minimizers to the assigned functions $\Psi_+$. This helps to derive a variety of uniform estimates for such solutions for all cones in $\mathcal{SC}_n$.

\begin{theorem}[Stable $L$-Vanishing on Cones]\label{miv} Let  $L$ be a natural and \si-adapted Schr\"odinger operator on cones $C \in   \mathcal{SC}_n$. Then, for any flat norm converging sequence $C_i \ra C_\infty$, $i \ra \infty$, with appropriately normalized associated solutions $\Psi_\pm(C_i)$ and $\Psi_\pm(C_\infty)$, we have
\[ \Psi_\pm(C_i) \circ \D \ra \Psi_\pm(C_\infty)\;\; C^{2,\alpha}\mm{-compactly on } C_\infty \setminus \sigma_{C_\infty} \mm{ as } i \ra \infty.\]
\end{theorem}

A more general  version of this result applies to ${\cal{H}}^{\R}_n$. A particularly interesting special case are degenerating sequences of regular Euclidean hypersurfaces $H_i \in {\cal{H}}^{\R}_n$ with singular limit. Then the induced solutions are always $L$-vanishing along the singular set of the limit hypersurface.

\begin{theorem}[Stable $L$-Vanishing on ${\cal{H}}^{\R}_n$]\label{miiv} Let $H_i \in {\cal{H}}^{\R}_n$ be a compactly converging sequence, $i\ge 1$, with limit $H_\infty \in {\cal{H}}^{\R}_n$. Further, let $L= -\Delta + V$ be a natural, \si-adapted Schr\"odinger operator with $\lambda^{\bp}_{L,H_i} \ge  c$ for some $c >0$ independent of $i$ and $u_i>0$ solutions of $L \, f =0$ on $H_i \setminus \Sigma_{H_i}$. Now assume the $u_i$ are $L$-vanishing along $U \cap \Sigma_{H_i}$, for some open $U\subset \R^{n+1}$. Then any induced solution $u>0$  on ${H_\infty}$ also $L$-vanishes along $U \cap \Sigma_{H_\infty}$.
\end{theorem}

The proofs of \ref{miii} - \ref{miiv} occupy the remainder of this section. The arguments for these results essentially coincide and will be addressed simultaneously. We first derive variants \ref{mi} of these results for minimal Green's functions which are $L$-vanishing along the entire singular set. This way we separate the invariance of  $L$-vanishing properties from the problem to work with the localization that the given solution $u >0$ of $L \, f =0$  is $L$-vanishing only in a neighborhood $U$ of some point $p \in \Sigma$. This localization is considered in a second step. This is another occasion where we employ the boundary Harnack inequalities already used to ensure the uniqueness of $\Psi_\pm$.\\

In these arguments we use some functional analytic consequences of \si-adaptedness. We show that the Hardy inequality extends from $C^\infty_0(H \setminus \Sigma)$ to $H^{1,2}_{\bp}(H \setminus \Sigma)$ and that furthermore, \si-adaptedness also allows us to control the $H^{1,2}_{\bp}$-norm, see \cite[Ch.5.1]{L1}.

\begin{lemma}[Equivalence of Norms]\label{hi2}
Let $L= -\Delta + V$ be a natural, \si-adapted Schr\"odinger operator with principal eigenvalue $\lambda^{\bp}_{L,H}>0$ for some $H \in \cal{G}$. Then there are constants $\beta^{\bp}_{L,H}$ and $\beta^{\bp,*}_{L,H} \ge 1$ such that
\begin{align}
\beta^{\bp}_{L,H} \cdot |f|_{H^{1,2}_{\bp}(H \setminus \Sigma)} & \ge \int_{H \setminus \Sigma}  f  \cdot  L f  \,  dV \, \ge \,  \lambda^{\bp}_{L,H}  \cdot \int_{H \setminus \Sigma} \bp^2\cdot f^2 dV\label{enq}\\
\beta^{\bp,*}_{L,H} \cdot \int_{H \setminus \Sigma}  f  \cdot  L f  \,  dV &\ge |f|_{H^{1,2}_{\bp}(H \setminus \Sigma)}\label{enq1}
\end{align}
for any $f \in H^{1,2}_{\bp}(H \setminus \Sigma)$.
\end{lemma}
\noindent\textbf{Proof} \,
We recall that $\bp$-adaptedness means that for some $a_L>0$ the potential $V$ satisfies
\begin{equation}\label{l}
 -a_L \cdot \bp^2 \le V  \le a_L \cdot \bp^2.
\end{equation}
Therefore,  \cite[Theorem 5.3]{L1}  shows that the inequality $\int_{H \setminus \Sigma}  f  \cdot  L f  \,  dV \, \ge \,  \lambda^{\bp}_{L,H}  \cdot \int_{H \setminus \Sigma} \bp^2\cdot f^2 dV$ not only holds for test functions in $C^\infty_0(H \setminus \Sigma)$ but actually for all functions in $H^{1,2}_{\bp}(H \setminus \Sigma)$. This implies the first inequality \eqref{enq} on $H^{1,2}_{\bp}(H \setminus \Sigma)$.\\

For the second inequality we use \eqref{l} again to write
\begin{align*}
(1 + (a_L+1)/\lambda^{\bp}_{L,H}) \cdot \int_{H \setminus \Sigma}  f  \cdot  L f  \,  dV  & \ge\\
\int_{H \setminus \Sigma} ( | \nabla_H f |^2 + V  \cdot  f^2 ) \, dV  \, + (a_L+1) \cdot \int_{H \setminus \Sigma} \bp^2\cdot f^2 dV & \ge  |f|_{H^{1,2}_{\bp}(H \setminus \Sigma)}
\end{align*}
whence \eqref{enq1} with $\beta^{\bp,*}_{L,H}:=1 + (a_L+1)/\lambda^{\bp}_{L,H}$. \qed

The subsequent theorem asserts that minimality of Green's functions $G(x,y)$ is stable under perturbation of the underlying area minimizer. In the argument we show that $G(x,y)$ can be described as finite energy minimizers of a Dirichlet type integral outside regular balls $B_\rho(x)$ centered around the pole $x$ of $G(x,\cdot)$. By the previous lemma this also entails finite $H^{1,2}_{\bp}$-norm. This is by no means obvious since in general, even minimal growth solutions strongly diverge towards $\Sigma$ and general  solutions usually lead to an infinite $H^{1,2}_{\bp}$-norm, cf.\ Theorem \ref{fix1} and Proposition \ref{evee} below for growth estimates.\\

We consider two situations of converging hypersurfaces with possibly non identical ambient spaces as described in Ch.\ref{tfr} and using $\D$-map identifications:
\begin{enumerate}[label=(\subscript{\textbf{S}}{{\textbf{\arabic*}}})]
\item  Let $H_i \in {\cal{H}}^{\R}_n$, $i\ge 1$, be compactly converging, with basepoints $a_i \in H_i \setminus \Sigma_{H_i}$ and possibly with $\Sigma_{H_i} \v$. Let $H_\infty \in {\cal{H}}$ be the limit with basepoint $a_\infty = \lim_{i \ra \infty} a_i$, $d(a_\infty ,0)=dist(a_\infty,\Sigma_{H_\infty}) = 5$ and $0 \in \Sigma_{H_\infty} \subset {\cal{H}}^{\R}_n$. Further, let $L= -\Delta + V$ be a natural, \si-adapted Schr\"odinger operator with $\lambda^{\bp}_{L,H_i} \ge  c$ for some $c >0$ independent of $i$.
\item  Let $H \in  {\cal{G}}^c$  and $H_i=\tau_i \cdot H \in {\cal{G}}^c$, for some sequence $\tau_i \ra \infty$. Also we choose basepoints $a_i \in  H \setminus \Sigma_H$, $i\ge 1$. We assume the $H_i$ and $a_i$ converge to some limit $H_\infty \in {\cal{H}}^{\R}_n$  with basepoint  $a_\infty = \lim_{i \ra \infty} a_i$, $d(a_\infty ,0)=dist(a_\infty,\Sigma_{H_\infty}) = 5$ and $0 \in \Sigma_{H_\infty} \subset {\cal{H}}^{\R}_n$. Further, let $L= -\Delta + V$ be a natural, \si-adapted Schr\"odinger operator on $H$.
\end{enumerate}

The reason why we only consider the case of blow-ups of one given almost minimizer $H \in {\cal{G}}^c$ is to ensure uniform control over the convergence and naturality properties of $\bp$.

\begin{proposition}[Stability of Minimal Green's Functions]\label{mi}
\noindent We assume we are in one of the two situations $S_1$ and $S_2$. We denote by $G_i=G_i(\cdot,a_i) >0$ the \textbf{minimal} Green's function on $H_i$ with pole in $a_i$ and normalized to $\int_{\p B_1(a_i)} G_i =1$, where $i \ge 1$ or $i=\infty$.\\

Then there is a subsequence, which we still denote $G_i$, such that as $i\ra\infty$,
\[
G_i \ra G_\infty\quad  C^{2,\alpha}\mm{-compactly on }H_\infty \setminus (\Sigma_{H_\infty} \cup \{a_\infty\})\mm{ via }\D\mm{-maps.}
\]
\end{proposition}

\noindent\textbf{Proof} \, We focus on case $S_1$.  Case $S_2$ then follows along similar lines noting that the scaling invariance of the $\bp^2$-weighted principal eigenvalue implies the condition
$\lambda^{\bp}_{L,H_i} \ge  c$ for some constant $c >0$ from the \si-adaptedness of $L$.\\

It is a trivial fact that the limit of the $G_i$ is again a Green's function; the point is to show its minimality. The idea is to characterize a minimal Green's function $G(\cdot,p)$  on $H$, outside some neighborhood $U$ of the pole $p \in H\setminus \Sigma_{H}$, as the minimizer of the following Dirichlet type integral
\[
J_{H \setminus U}(f):=  \int_{H \setminus U}   | \nabla_H f |^2 + V  \cdot  f^2 \, d V,\;f \mm{ smooth with }supp \, f \subset H\setminus \Sigma_{H}\mm{ and }f|_{\p U} = G.
\]
\si-adaptedness of $L$ implies then that the minimality of the $J$-integral is preserved under convergence of the underlying hypersurfaces.\\

In our situation where we have a compactly converging sequence of pointed minimizers $H_i \ra H_\infty$ with $a_i \ra a_\infty$ we set $U_i=\D(B_1(a_\infty))$ and $U_\infty= B_1(a_\infty)$.  Without loss of generality we may also assume that $B_{5}(a_i) \subset H_i \setminus \Sigma_{H_i}$, since  $B_{5}(a_\infty) \subset H_i \setminus \Sigma_{H_\infty}$, and these balls are, possibly after rescaling, diffeomorphic and uniformly almost isometric in $C^3$-norm to the Euclidean ball $B_{5}(0)$.\\

The proof will be divided into three steps. In the first two preliminary steps we derive results needed in the third and main step. In the first step we show that the functionals $J_{H_i\setminus U_i}$, $i\geq0$ or $i=\infty$, are uniformly bounded. In Step 2 we shall prove that the $G_i$ are the unique minimizers of the Dirichlet type integral above. In Step 3 we prove that $G_i$ converge to a limit solution $G^*$ which minimizes the Dirichlet integral on $H_\infty \setminus U_\infty$. (The uniqueness from Step 2 then implies that $G^*= G_\infty$.) For this we prove that there is no bubbling-off of negative contributions to $J_{H_i \setminus U_i}(f)$ before we reach the limit. This will use the \si-adaptedness of $L$ and the fact  that
$H^{1,2}_{\bp}(H \setminus \Sigma) \equiv H^{1,2}_{\bp,0}(H \setminus \Sigma)$ proved in \cite[Theorem 5.3]{L1}.\\

\textbf{Step 1} \,
\emph{For the infima $\inf J_i$ of $J_{H_i\setminus U_i}(\cdot)$, $i\geq 0$ or $i=\infty$, on $H^{1,2}_{\bp}(H_i\setminus \Sigma_{H_i})$, we have a common constant $\alpha>0$ so that}
\[|\inf J_i| \le \alpha, \mm{ for any } i\geq 0 \mm{ or } i=\infty.\]

For the \emph{lower} bound let $f$ be a smooth function with $supp \,f \subset H_i\setminus \Sigma_{H_i}$ and $f|_{\p U_i} = G_i$. We extend $f$ as a smooth function to $F_i$ on $U_i$. Hence $|J_{U_i}(F_i)|  \le c$ for some uniform constant $c < \infty$. The Hardy inequality for $L$ entails, by Lemma \ref{inhnat}, that $J_{H_i}(F_i) \ge 0$, whence $J_{H_i\setminus U_i}(f) \ge J_{H_i}(F_i) -c \ge -c$. On the other hand, a fixed test function $f$ supported in  $B_2(a_\infty)$ gives a uniform \emph{upper} bound for $J_{H_\infty \setminus U_\infty}(f)$ on $H_\infty$ and thus, via almost isometric $\D$-map pull-back for sufficiently large $i$, also for $J_{H_i \setminus U_i}(f)$ on $H_i$.\qed\\

\textbf{Step 2} \,
\emph{For the Dirichlet problem of $J_{H_i\setminus U_i}(\cdot)$ over $H^{1,2}_{\bp}(H_i\setminus \Sigma_{H_i})$ the minimal Green's function $G_i$ restricted to $H_i\setminus U_i$ is the unique minimizer $v>0$ with boundary value $\phi= G_\infty \circ \D$ on $\p U_i$.}

\begin{remark}\label{tra}1. If we work with functions in $H^{1,2}_{\bp}$, then the Dirichlet problem is considered in the \emph{trace sense}, cf.\ \cite[Chapter 5.5 and 6.1]{E}. At any rate, the resulting minimizers are regular and thus solve the problem in the classical sense. \\

2. The claim follows from the general theory of symmetric semi-bounded operators cf. Remark \ref{frr} below. But the following explicit argument may be more satisfactory.\qed
\end{remark}

For the uniqueness part, we assume we had two distinct positive solutions $v_1 \neq v_2$ with $J_{H_i\setminus U_i}(v_1)=J_{H_i \setminus U_i}(v_2)= \inf J_{H_i\setminus U_i}$. Regularity theory shows that $v_{1,2}$ are at least $C^{2,\alpha}$-regular. Upon relabeling $v_1$ and $v_2$ we can find an $\eta <1$ such that $\eta \cdot v_1$ and $v_2$ coincide in a non-empty set $C$ outside $\overline{B_1(p)}$. The Hopf maximum principle shows that $C$ must be a smooth submanifold where the graphs of $v_1$ and $v_2$ intersect transversally. Since both functions minimize $J$ we observe that
\[
J_{H_i\setminus U_i}(\min\{\eta  \cdot v_1,v_2\})= \eta^2 \cdot \inf J_{H_i\setminus U_i} \, \mm{  or } \,  J_{H_i\setminus U_i}(\max\{\eta \cdot v_1,v_2\})= \inf J_{H_i\setminus U_i}.
\]
Thus $\min\{\eta  \cdot v_1,v_2\}$ or $\max\{\eta \cdot v_1,v_2\}$  also minimizes the functional $J$ and must be smooth. But this contradicts the fact that both $\min\{\eta  \cdot v_1,v_2\}$ and $\max\{\eta \cdot v_1,v_2\}$ are non-smooth along $C$. This proves the uniqueness assertion.\\

For the existence part we choose a sequence of smoothly and compactly bounded domains $\overline{D}_m\subset D_{m+1} \subset H_i\setminus \Sigma_{H_i}$, $m \ge 0$, with $\bigcup_m D_m = H_i\setminus \Sigma_{H_i}$ and $\overline{U} \subset D_1$. There is a unique solution  $Q_{i,m}$ for the Dirichlet problem with boundary data $Q_{i,m} \equiv G_i$ on $\p U$ and $Q_{i,m} \equiv 0$ on $\p D_m$ which minimizes the Dirichlet functional. We infer from \cite[Theorem 5.3]{L1} and the proof of the critical operator case in \cite[Theorem 5.7]{L1} that $Q_{i,m} \nearrow G_i$ and $J_{H_i\setminus U_i}(Q_{i,m}) \searrow  \inf J_{H_i\setminus U_i}$ as well as $J_{H_i\setminus U_i}(Q_{i,m}) \searrow J_{H_i\setminus U_i}(G_i)$. Thus the unique minimizer of $J_{H_i\setminus U_i}$ is the minimal Green's function. \qed\\

For any $\eta >0$, comparing the functional over compactly supported test functions on the limit space $H_\infty$ with the functional over compactly supported test functions on $H_i$ via $\D$-maps yields an $i_\eta$ with
\begin{equation}\label{e}
 J_{H_i \setminus U_i}(G_i) \le \eta+  J_{H_\infty \setminus U_\infty}(G_\infty)\mm{ for any } i \ge i_\eta.
\end{equation}
This, in turn, implies that
\begin{equation}\label{e1}
|G_i|_{H^{1,2}_{\bp}(H_i \setminus (\Sigma_{H_i} \cup  U_i))} \le \beta^{\bp,*}_{L,H} \cdot (\eta +  J_{H_\infty \setminus U_\infty}(G_\infty)+c^*),
\end{equation}
cf.\ \eqref{enq1} of Lemma \ref{hi2}.  $c^*>0$ is a common upper bound for the contribution of the extensions of the $G_i$ to $U_i$ by bounded regular functions. By regularity theory the $G_i$ compactly $C^{2,\alpha}$-subconverge to some solution $G^*>0$ with  $G^* \in H^{1,2}_{\bp}(H_\infty \setminus (U_\infty \cup \Sigma_{H_\infty}))$ and $G^* = G_\infty \mm{ on } \p U$.\\

\textbf{Step 3} \,
\emph{$G^*$ minimizes the Dirichlet integral whence, $G^*=G_\infty$.}\\

In view of the compactly $C^{2,\alpha}$-subconverging $G_i$, this amounts to exclude bubbling-off phenomena of negative contributions to the Dirichlet integral concentrating near $\Sigma_i$ which are no longer visible in the limit. Towards this end, we let $\I_\ve:= \{x \in H_i \setminus \Sigma_{H_i} \, |\, \delta_{\bp}(x) < \ve\}$, for $\ve>0$ which we choose sufficiently small so that $U_i \cap \I_{3 \cdot \ve} \v$ and $R \ge 10$. Assume that there is a constant $c>0$ independent of both, $\ve$ and $R$, such that there is an $i^0_{\ve.R}$ so that
\[
J_{H_i\setminus U_i}(G_i) = J_{B_{2 \cdot R}(a_i)\setminus (U_i \cup \I_\ve)}(G_i) + J_{(H_i \setminus B_{2 \cdot R}(a_i)) \cup \I_\ve}(G_i) \mm{ and }  J_{(H_i \setminus B_{2 \cdot R}(a_i)) \cup \I_\ve}(G_i) \le - c
\]
for any $i \ge i^0_{\ve,R}$. To argue that this cannot happen we use the setup from the proof of \cite[Theorem 5.3]{L1}. We first construct a cut-off function concentrated near $\Sigma_i$ by taking some fixed $\psi \in C^\infty(\R,[0,1])$ with $\psi \equiv 1$ on $\R^{\le 0}$ and $\psi \equiv 0$ on $\R^{\ge 1}$. Then we set
\begin{equation}\label{cut}
\psi[\zeta](x):= \psi(\zeta^{-1} \cdot \delta_{\bp}(x)-1)\mm{ for } x \in H_i \setminus \Sigma_i \mm{ and } \zeta \in (0,1).
\end{equation}
Since $\delta_{\bp}$ is Lipschitz continuous there exists some constant $c(\psi) >0$ depending only on $\psi$ such that for the distributional derivative
\begin{equation}\label{nabb}
|\nabla \psi[\zeta](x)| \quad
\begin{cases}
\le c(\psi) \cdot \bp(x), &\text{for $\delta_{\bp}(x) \in (\zeta, 2 \cdot \zeta)$;}\\
= 0, &\text{otherwise.}
\end{cases}
\end{equation}
From \eqref{cut} and \eqref{nabb} we see  that
\begin{align*}
|\psi[\zeta] \cdot G^*|_{H^{1,2}_{\bp}(\I_{2\zeta} \setminus \I_\zeta)}=&\int_{\I_{2\zeta} \setminus \I_\zeta}|\nabla (\psi[\zeta] \cdot G^*)|^2 + \bp^2 \cdot (\psi[\zeta] \cdot G^*)^2 \, dV\\
\le&\int_{\I_{2\zeta} \setminus \I_\zeta} 2|\nabla \psi[\zeta]|^2  \cdot (G^*)^2 + 2|\nabla G^*|^2 + 2\bp^2  \cdot (G^*)^2 \, dV\\
\le&\int_{\I_{2\zeta} \setminus \I_\zeta} 4c(\psi)^2 \cdot \bp^2  \cdot (G^*)^2 + 2|\nabla G^*|^2 + 2\bp^2  \cdot (G^*)^2 \, dV\\
\le& \, 4(c(\psi)^2 +1) \cdot   |G^*|_{H^{1,2}_{\bp}(\I_{2 \cdot \zeta} \setminus \I_\zeta)}
\end{align*}

Secondly, we consider a cut-off towards infinity by setting $\psi_R(x):= \psi(R^{-1} \cdot |x|-1)$ for $x \in H_\infty \setminus \Sigma_\infty$ and $R \ge 10$. We note as a counterpart of \eqref{nabb}  that
\begin{equation}\label{nabba}
|\nabla \psi_R(x)| \quad
\begin{cases}
\le c(\psi)^*/R, &\text{for $|x| \in (R, 2 \cdot R)$;}\\
= 0, &\text{otherwise}
\end{cases}
\end{equation}
for some $c(\psi)^* >0$. The Lipschitz condition $  |\delta_{\bp_{H_\infty}}(x)- \delta_{\bp_{H_\infty}}(a_\infty)|   \le L_{\bp} \cdot d_{H_\infty}(x,a_\infty) $  shows that $\delta_{\bp_{H_\infty}}(x)  \le L_{\bp} \cdot d_{H_\infty}(x,a_\infty) + \delta_{\bp_{H_\infty}}(a_\infty)$. This means  that
\begin{equation}\label{ab}
\bp(x) \ge b_{H_\infty} \cdot d_{H_\infty}(x,a_\infty) ^{-1} \mm{ on } H_\infty \setminus (\Sigma_\infty \cup B_1(a_\infty)),
\end{equation}
for some $b_{H_\infty}>0$. Then we have using  \eqref{nabba} and \eqref{ab}:

\begin{align*}
|(\psi_R \cdot G^*|^2_{H^{1,2}_{\bp}(B_{2 \cdot R} \setminus B_R(a_\infty))}&= \int_{B_{2 \cdot R} \setminus B_R(a_\infty)}|\nabla (\psi_R \cdot G^*)|^2 + \bp^2 \cdot (\psi_R \cdot G^*)^2 \, dV\\
&\le 2  \int_{B_{2 \cdot R} \setminus B_R(a_\infty)}\psi_R^2 \cdot | \nabla G^*|^2 + |\nabla \psi_R|^2\cdot (G^*)^2 + \bp^2 \cdot \psi_R^2 \cdot (G^*)^2 \, dV\\
&\le 2 \int_{B_{2 \cdot R} \setminus B_R(a_\infty)}\psi_R^2 \cdot | \nabla G^*|^2 + \bp^2 \cdot \psi_R^2 \cdot(G^*)^2 \, dV \\
&\quad  +  2  \int_{B_{2 \cdot R} \setminus B_R(a_\infty))}  (c(\psi)^*/R)^2 \cdot (G^*)^2 \, dV\\
&\le    2 \cdot \left(1 +(c(\psi)^*/b_{H_\infty})^2\right) \cdot |G^*|_{H^{1,2}_{\bp}(B_{2 \cdot R} \setminus B_R(a_\infty))}.
\end{align*}

We observe, from $G^* \in H^{1,2}_{\bp}(H_\infty \setminus (U_\infty \cup \Sigma_{H_\infty}))$, that  $|G^*|_{H^{1,2}_{\bp}(\I_{2 \cdot \zeta} \setminus \I_\zeta)}$
and $|G^*|_{H^{1,2}_{\bp}(B_{2 \cdot R} \setminus B_R(a_\infty))}$   tend to $0$, when $\ve \ra 0$ respectively   $R \ra \infty$. Transferring the cut-off region \[C[R,\ve]:=\big(B_{2 \cdot R}(a_\infty) \cap \I_{2 \cdot \ve} \setminus \I_\ve \big) \cup \big(B_{2 \cdot R} \setminus B_R(a_\infty)\big)  \setminus \I_\ve\] via the $\D$-map to $H_i$ shows that for sufficiently small $\ve>0$, large $R$ and $i$, the contributions to the two norms on the cut-off region become arbitrarily small:
\[
|\psi[\ve] \cdot G_i|_{H^{1,2}_{\bp}(\D (C[R,\ve]))} \ll 1 \mm{ and thus }
|J_{(\D (C[R,\ve]))}(\psi[\ve] \cdot G_i)| \ll 1
\]
due to the $\bp$-adaptedness of $L$. That is, we get for large $i$
\[
(1-\psi_R) \cdot \psi[\ve] \cdot G_i \in H^{1,2}_{\bp}(H_i \setminus \Sigma_{H_i})\quad\mm{and}\quad J_{H_i}((1-\psi_R) \cdot \psi[\ve] \cdot G_i) < -c/2.
\]
But this contradicts the \si-adaptedness of $L$ for test functions in $H^{1,2}_{\bp}(H_i \setminus \Sigma_{H_i})$ as this asserts that $J_{H_i}((1-\psi_R) \cdot \psi[\ve] \cdot G_i)>0$, cf.\ Lemma \ref{hi2}. \qed

\begin{remark}\label{frr}
The considered domain of a shifted \si-adapted operator $L$ is typically a space of (sufficiently) regular functions. For such domains the operators are symmetric on the original space $H \setminus \Sigma$, but in general \emph{not} self-adjoint. (Note in passing that the associated operator $\delta_{\bp}^{-2} \cdot L$ on the hyperbolic unfolding is usually no longer a symmetric operator.)\\

However, right from the definition, any shifted \si-adapted operator $L$ is a symmetric \emph{semi-bounded} operator on $C^\infty_0(H \setminus \Sigma)$. This is sufficient to ensure that $L$ admits a canonical \emph{self-adjoint extension} $L_F$, the so-called \emph{Friedrichs extension}. Its domain $D(L_F)$ satisfies
\[
C^\infty_0(H \setminus \Sigma) \subset D(L_F) \subset H^{1,2}_{\bp}(H \setminus \Sigma).
\]
For these Friedrichs extensions we have a general existence and uniqueness theory for the Dirichlet problem in Step 2 of Proposition~\ref{mi} above, cf.\ \cite[Ch.4.3]{Hf} and \cite[Ch.5.5]{Z}. \qed
\end{remark}

The boundary Harnack inequality allows us to derive the following localized versions of this stability of minimal Green's functions.

\begin{proposition}[Localized Inheritance of $L$-Vanishing]\label{mj}
\noindent We assume situation $S_1$ or $S_2$. Then, let $u_i >0$ be a solution of $L\,f=0$ on each $H_i$ which via $\D$-maps are compactly converging to some entire solution $u_\infty>0$ of $L\,f=0$ on $H_\infty$. If there is a ball $B$ in a common ambient manifold of the $H_i$ and $H_\infty$ around a singular point $p \in \Sigma_{H_\infty}$ such that $u_i$, for any $i$, is $L$-vanishing along $B \cap \Sigma_{H_i}$, then $u_\infty$ is also $L$-vanishing along $B \cap \Sigma_{H_\infty}$.
\end{proposition}

The proof is done in three steps. In Step 1 we prepare canonical $\Phi_\delta$-chains on the $H_i$ converging to a canonical $\Phi_\delta$-chains on $H_\infty$ we use as a neighborhood basis of $0$. Next, in Step 2, we observe that for these  $\Phi_\delta$-chains on the $H_i$ we get common constants in the BHP on all $H_i$, cf.\cite[Ch.2.1]{L1}.  In Step 3 we use the stability of minimal Green's functions $G_i$ on the $H_i$ \ref{mi}. Using Step 2, we see that the minimal Green's function $G_i$ upper bounds the $u_i$ on $B \cap H_i$ by a common multiple of $G_i$ and this carries over to the limit. Thus, we infer that  $u_\infty$ is also $L$-vanishing along $B \cap \Sigma_{H_\infty}$.

\begin{remark}\label{qia} For Step 1 we recall from \cite[Ch.3.1]{L1} that the BHP for \si-adapted operators is a direct consequence of Ancona's BHP from the unfolding correspondence. For converging sequences $H_i$ as in situations $S_1$ and $S_2$ the associated Gromov hyperbolic  $(H_i,d_{\bp_{H_i}})$
converge compactly to that of the limit space $H_\infty$. For the Whitney smoothed \si-metrics  $d_{\bp^*}$, used in the unfolding correspondence, any such sequence $(H_i,d_{\bp^*_{H_i}})$ has subsequences converging to a Whitney smoothing $(H_\infty,d_{\bp^*_{H_\infty}})$. This readily follows from the proof of \cite[ Proposition B.3]{L1} since the basepoints of the \si-adpated covers compactly subconverge. In what follows we therefore assume that the $(H_i,d_{\bp^*_{H_i}})$ already converge to $(H_\infty,d_{\bp^*_{H_\infty}})$.\qed
\end{remark}

\noindent\textbf{Proof of \ref{mj}} \,  \textbf{Step 1} \, We first recall the construction of \emph{canonical $\Phi_\delta$-chains} from  \cite[Lemma 2.4]{L1}. For a $\delta$-hyperbolic space we consider geodesic arcs $\gamma: (0,  c)  \ra  X$, for $c>10^3 \cdot \delta$. We choose the track points $x_k=\gamma(t_k)$, $k=0,...,m+1$, with $x_0=\gamma(0), x_{m+1}=\gamma(c)$, with $d(x_k,x_{k+1})=300\cdot \delta$, for $k <m$, and $d(x_m,x_{m+1})\le 300 \cdot \delta$. Then the $\mathcal{N}^\delta_k(\gamma):=U_{t_k}$ with
\begin{equation}\label{pch}
U^\gamma_t:=\{x \in X\,|\, dist\big(x,\gamma([t,c))\big)<dist\big(x,\gamma((0,t]) \big)\}.
\end{equation}
 form a canonical $\Phi_\delta$-chain with track points $x_k ,k=1,...,m$,
for $\Phi_\delta(t)=a_\delta + b_\delta \cdot t$ for some $a_\delta, b_\delta  > 0$ depending only on $\delta$. We call $m$ its length.\\

In both situations $S_1$ and $S_2$ we now turn to the Whitney smoothed hyperbolic unfoldings. We can choose hyperbolic geodesic arcs $\gamma[c_i]:(0,c_i)\ra  H_i \setminus \Sigma_{H_i}$, with basepoint $\gamma[c_i](0)= c_i$ and $c_i \ra \infty$ for $i \ra \infty$, so that the $\gamma[c_i]$ compactly converge to some $\gamma[c_\infty]  \subset (H_\infty,d_{\bp^*_{H_\infty}})$ representing the Gromov boundary point $0 \in \Sigma_{H_\infty}$. The convergence is formalized via $\D$-maps.\\
This also implies that the canonical $\Phi_\delta$-chains $\mathcal{N}^\delta_k(\gamma[a_i])\subset H_i$, of length $m(i) \ra \infty$, we assign to each of the  $\gamma[a_i]$, compactly converge to  $\mathcal{N}^\delta_k(\gamma[a_\infty])\subset H_\infty$, via $\D$-maps.   \\

\textbf{Step 2} \, We can now apply the BHPs \cite[Theorem 3.4 and 3.5]{L1}. For a minimal Green's function $G$ and solutions $v >0$ of $L\, f = 0$ on $H_i \setminus \Sigma_{H_i}$ which is  \emph{L}-vanishing along $\mathcal{N}^\delta_2(\gamma[a_i])$ we have:
\[G(x,p)/v(x) \le C \cdot  G(y,p)/v(y),\mm{ for all } p \in (H_i \setminus \Sigma_{H_i}) \setminus \mathcal{N}^\delta_1(\gamma[a_i]) \mm{ and } x,\, y \in \mathcal{N}^\delta_3(\gamma[a_i]).\]
For $k_L \le \kappa$ and $\ve_L \ge \eta$ the constant $C$ only depends on $\kappa$ and $\eta> 0$. In our cases  this means that we can choose a common Harnack constant $C$ for all $i$. For $S_1$ this from the stable BHP on ${\cal{H}}^{\R}_n$ \cite[Theorem 3.5]{L1}. It gives the same constant for all such hypersurfaces. For $S_2$ we can employ the BHP on ${\cal{G}}_n$ of \cite[Theorem 3.4]{L1} since we are working only with one fixed hypersurface.\\

\textbf{Step 3} \, In both situations, $S_1$ or $S_2$, we have $u_i \le C \cdot G_i$   on  $\mathcal{N}^\delta_3(\gamma[a_i])$ for some common $C\ge 1$.
Now we recall from \cite[Proposition 8.10]{BHK} that the canonical $\Phi_\delta$-chain $\mathcal{N}^\delta_k(\gamma[a_\infty])\subset H_\infty$ for $\gamma[a_\infty]$ forms a neighborhood basis of $0 \in H_\infty$. Since $\Sigma_{H_\infty} \cap \overline{B}$ is compact, we can therefore assume (without loss of generality) that $B \cap  H_\infty \subset  \mathcal{N}^\delta_k(\gamma[a_\infty])$,
for $k=1,...3$,  $i \ge 1$ and $i=\infty$.\\

Thus we also have
\[u_\infty \le C \cdot G_\infty \mm{ on } B \cap H_\infty \setminus \Sigma_{H_\infty}\]
Since $G_\infty$ is minimal we see that $u_\infty$ also $L$-vanishes along $B \cap \Sigma_{H_\infty}$. \qed

This result also concludes the \textbf{proof} of \ref{miii} and \ref{miv} taking increasingly larger balls $B$ in \ref{mj}.

%
%
%
\setcounter{section}{4}
\renewcommand{\thesubsection}{\thesection}
\subsection{Martin Theory on Minimal Cones}

Here we refine our analysis of natural Schr\"odinger operators. The minimal growth stability and Martin theory can be used to build an inductive asymptotic analysis of solutions with minimal growth near singular points using iterative blow-ups. In the cases of the Jacobi field operator and the conformal Laplacian we derive some further details.

\subsubsection{Schr\"odinger Operators and Scaling Actions} \label{nee}
For Schr\"odinger operators $L$ which are naturally associated with $H\in \cal{G}$ we get a neat representation for certain distinguished solutions of $L\,f=0$. We estimate their radial growth in the  cone case to gain a detailed understanding of the growth of the Martin kernel for compact $H$.

\begin{definition}[Natural Schr\"odinger Operators]\label{spli}
A natural and shifted \si-adapted operator $L$ is called a \textbf{natural Schr\"odinger operator} if for any given $H \in \cal{H}$ the operator $L(H)$ has the form
\[
L(H)(u)= -\Delta_H \, u + V_H(x) \cdot \, u \mm{ on } H \setminus \Sigma_H
\]
for some H\"older continuous function $V_H(x)$.
\end{definition}

\begin{remark}[Inheritance versus Naturality]\label{ivn}
Note that a natural Schr\"odinger operator is merely shifted \si-adapted with respect to a given $H \in \cal{H}$. However, due to Inheritance Lemma~\ref{inhnat}, naturality always implies that for a blow-up $N$ of $H$ around a singular point in $\Sigma_H$, $L(N)$ is shifted \si-adapted with $\lambda^{\bp}_{L,H} \le \lambda^{\bp}_{L,N}$.
\qed
\end{remark}

This class of Schr\"odinger operators  covers the examples we discussed in Chapter \ref{exp}, namely the Laplacian $-\Delta_H$, the Jacobi field operator $J_H=-\Delta_H - |A|^2-Ric_M(\nu,\nu)$ and the $S$-conformal Laplacian $L_{H,S} = -\Delta_H +\frac{n-2}{4 (n-1)}\cdot (scal_H - S|_H)$.\\

For an area minimizing cone $C \subset \cal{H}$ the naturality implies that $V_C(t \cdot x)=t^{-2} \cdot V_C(x)$ for any $x \in C \setminus \sigma_C$ and $t >0$.
In this case we can write $V(x)=r^{-2} \cdot V^\times(\omega)$ for $x=(\omega,r) \in S_C \setminus \sigma_C \times \R^{> 0} = C \setminus \sigma_C$. It follows that
\begin{equation}\label{pol}
L \, v= -\frac{\p^2 v}{\p r^2} - \frac{n-1}{r} \cdot \frac{\p v}{\p r} -
\frac{1}{r^2} \cdot \Delta_{S_C} v  +\frac{1}{r^2} \cdot  V^\times(\omega) \cdot v=: -\frac{\p^2 v}{\p r^2} - \frac{n-1}{r} \cdot \frac{\p v}{\p r} +\frac{1}{r^2} \cdot L^\times \, v,
\end{equation}
where $L^\times=L^\times(S_C)$ is an operator on $S_C$. We will see in \ref{fix1} below that $L^\times$ is also a natural Schr\"odinger operator.\\

\noindent\textbf{Scaling Actions} \,
On an area minimizing cones equipped with a natural and \si-adapted  operator $L$, both, the operator and solutions of $L \, f =0$ reproduce themselves under scalings of the cone up to constant multiples. Concretely, we express $L$ on $C = \R^{\ge 0} \times S_C$ in geodesic coordinates $x_1=r$, $x_2,...,x_n$ such that $x_2,...,x_n$ locally parametrize $S_C=C\cap\p B_1(0)$:
\[
-L(u) =  \sum_{i,j} a_{ij} \cdot \frac{\p^2 u}{\p s_i \p s_j} + \sum_i b_i
\cdot \frac{\p u}{\p s_i} + c \cdot u.
\]
When $L$ is natural this means that for any $\eta>0$,
\[
a_{ij}(\eta \cdot x) = a_{ij}(x),\, b_i(\eta \cdot x) = \eta^{-1} \cdot b_i(x) \mm{ and } c(\eta \cdot x) = \eta^{-2} \cdot c(x).
\]
Thus, for any function $u(x)$ that solves $L \, f =0$ its rescaling $u(\eta \cdot x)$ also solves this equation by the chain rule. In particular, the Green's function and the set of minimal solutions of $L \, f = 0$ are reproduced up to multiples under composition with the scaling map
\[
S_\eta: C \ra C \mm{ given by } x \mapsto \eta \cdot x \mm{ for } \eta \in (0,\infty).
\]
More concretely, we consider the map $u \mapsto u \circ S_\eta$ and normalize the values of the resulting functions $u \circ S_\eta$ to $1$ in some base point $p \in C \setminus \sigma_C$. In this way we define a scaling action $S^*_\eta$ on the Martin boundary.

\begin{lemma}[Attractors and Fixed Points of $\mathbf{S^*_\eta}$]\label{fix}
Let $C \in \mathcal{SC}_n$ be a singular area minimizing cone and $L$ be a natural \si-adapted operator on $C$. Then we have:
\begin{enumerate}
  \item The scaling action $S^*_\eta$ on $\p_M (C,L)$ has exactly two fixed points:  the tip $0_C$ and the point at infinity $\infty_C$, both viewed as minimal functions.
  \item For $z \in \p_M (C,L) \setminus \{0_C,\infty_C\}$ we find
\[
S^*_\eta(z) \ra \infty \mm{ as } \eta \ra \infty_C \, \mm{ and } \,  S^*_\eta(z) \ra 0_C \mm{ as } \eta \ra 0.
\]
  \item More generally, for any given $z \in C \setminus \sigma_C$ we have
  \[
  \frac{G(x,S_\eta(z))}{G(p,S_\eta(z))} \ra \infty_C \mm{ as } \eta \ra \infty \,\mm{ and }\, \frac{G(x,S_\eta(z))}{G(p,S_\eta(z))} \ra 0_C \mm{ as } \eta \ra 0
  \]
  locally uniformly for $x \in C \setminus \sigma_C$.
\end{enumerate}
\end{lemma}
\noindent\textbf{Proof} \,
This readily follows from the way the pole of $G(\cdot,S_\eta(z))$ shifts under these scaling operations. \qed\\

Now we give a description of the two fixed point solutions in $\p_M (C,L)$. Despite the
fact mentioned above that there is no proper way of transforming the Martin theory from the original hypersurface $H$ to its tangent cones this will allow us to build an inductive decomposition scheme for solutions on $H$.

\begin{theorem}[Separation of Variables]\label{fix1}
For any cone $C \in \mathcal{SC}_n$ and any natural Schr\"odinger operator $L$ let us consider the \si-adapted operator $L_\lambda=L - \lambda \cdot \bp^2 \cdot Id$,  $\lambda < \lambda^{\bp}_{L,C}$. Then we have:
\begin{itemize}
  \item Viewed as functions $\Psi_-=0_C$ and $\Psi_+=\infty_C$ on $C \setminus \sigma$ the two fixed points $0_C$, $\infty_C \in \p_M (C,L_\lambda)$ can be written as
\[
\Psi_\pm(\omega,r) = \psi(\omega) \cdot r^{\alpha_\pm}
\]
with $(\omega,r) \in S_C \setminus \sigma \times \R$ and $\alpha_\pm = - \tfrac{n-2}{2} \pm \sqrt{ \Big( \tfrac{n-2}{2} \Big)^2 + \mu_{C,L^\times_\lambda}}$.
  \item Let $\bp^\times(\omega)=r \cdot \bp(x)$ for $x=(\omega,r) \in C \setminus \sigma$. The associated operator
\[
L_\lambda^\times(v)(\omega)= - \Delta_{S_C} v(\omega) + \big(V^\times(\omega) - \lambda \cdot\bp^\times(\omega)^2 \big)\cdot v(\omega)
\]
defined on $S_C \setminus \Sigma_{S_C}$ is a natural Schr\"odinger operator with non-weighted principal eigenvalue $\mu_{C,L^\times_\lambda} > - (\frac{n-2}{2})^2$ and ground state $\psi(\omega) >0$, that is,
\[
L_\lambda^\times \, \psi =  \mu_{C,L^\times_\lambda} \cdot  \psi, \mm{ on } S_C \setminus \Sigma_{S_C}.
\]
\end{itemize}
\end{theorem}
\noindent\textbf{Proof} \,
We proceed in three steps: First, we show that $\Psi_\pm$ can be written as a product $\Psi_\pm(\omega,r) = \psi_\pm(\omega) \cdot r^{\alpha_\pm}$. Then we determine $\psi_\pm$ and $\alpha_\pm$ for some inner approximation $C$ by regular subcones. Finally, we prove that the resulting values converge to $\psi_\pm$ and $\alpha_\pm$ on $C$.\\

\noindent\textbf{Product Shape} \,
We first restrict the two fixed points $\Psi_\pm \in \p_M (C,L_\lambda)$ to a regular ray $\Gamma_v = \R^{>0} \cdot v \in C \setminus \sigma$, for some $v\in S_C$. Consider the map $\Psi_\pm|_{\Gamma_v}: \R^{>0} \cdot v \ra \R^{>0}$ which we view as a restriction of a function $f:\R^{>0} \ra \R^{>0}$. Now up to a constant $\Psi_\pm$ reproduces under scalings: For any $\eta >0$, there is constant $c_\eta>0$ such that
\[
f(\eta \cdot x)=c_\eta \cdot f(x) \mm{ for all } x \in \R^{>0}.
\]
From this it follows that $f$ is a monomial, that is, $f(x)=a \cdot x^b$ for some constants $a(v) > 0$, $b(v) \in\R$. This argument applies to any regular ray $\Gamma_v$.\\

Next we consider the Harnack inequality for $L_\lambda$ on a ball $B_{2R}(v) \subset C \setminus \sigma_C$ for some $R>0$. We get, for any positive solution $u$ of
$L_\lambda \, f =0$, the Harnack inequality
\begin{equation}\label{h1}
\sup_{B_R(v)}\, u \le c(L_\lambda,v,R)  \cdot \inf_{B_R(v)}\, u,
\end{equation}
for some constant independent of $u$. The crux of the matter is that both the scaling symmetry of $C$ and the naturality of $L_\lambda$ imply that the same constant $c$ can still be used in the Harnack inequality after scalings around the tip $0$: For any $s_0$ we have
\begin{equation}\label{h2}
\sup_{B_{s \cdot R}(s \cdot v)}  u \le c(L_\lambda,v,R)  \cdot \inf_{B_{s \cdot R}(s \cdot v)} u.
\end{equation}
This implies that for all rays $\Gamma_w$ passing through $B_R(v)$ the exponent $b(w)$ equals $b(v)$. Since $S_C \setminus \Sigma_{S_C}$ is connected, this shows that $b(v)$ is constant on $S_C \setminus \Sigma_{S_C}$. Thus $\Psi_\pm$ can be written as $\Psi_\pm(\omega,r) = \psi_\pm(\omega) \cdot r^{\alpha_\pm}$ for $(\omega,r) \in \p B_1(0) \cap C \setminus
\sigma \times \R^{> 0}$ and some  $C^{2,\beta}$-regular function $\psi_\pm$ on $\p B_1(0) \cap C \setminus \sigma_C \times \R^{> 0}$, $\beta>0$.\\

We reinsert $\Psi_\pm(\omega,r) = \psi_\pm(\omega) \cdot r^{\alpha_\pm}$ into the equation $L_\lambda \, f =0$ which we write in polar coordinates as in \eqref{pol}. A separation of variables shows that the $\psi_\pm$ solve the equations
\begin{equation}\label{sph}  L_\lambda^\times \,
v := - \Delta_{S_C} v(\omega) + \big(V^\times(\omega) - \lambda \cdot (\bp^\times)^2(\omega) \big)\cdot v(\omega) = (\alpha_\pm^2 + (n-2)\cdot  \alpha_\pm) \cdot v
\end{equation}
on  $S_C \setminus \Sigma_{S_C}$. Further, $L_\lambda^\times$ is again a natural Schr\"odinger operator and adapted to the \si-transform $\bp^\times=\bp_C|_{S_C}$.\\

\noindent\textbf{Inner Regular Approximation} \,
We use an approximation by Dirichlet eigenvalue problems to show that
\begin{itemize}
  \item $\alpha_-^2 + (n-2)\cdot  \alpha_- = \alpha_+^2 + (n-2)\cdot  \alpha_+ > - (n-2)^2/4$.
  \item $\mu_{C,L^\times_\lambda} := \alpha_\pm^2 + (n-2)\cdot  \alpha_\pm $ is the non-weighted  \emph{principal eigenvalue} of $L_\lambda^\times$. Further, $L_\lambda^\times$ is shifted \si-adapted.
  \item $\psi_-=\psi_+$  and $\psi:=\psi_\pm$ is the corresponding \emph{ground state} of $L_\lambda^\times$.
\end{itemize}
Towards this end we choose smoothly bounded domains $D_i \subset \overline{D}_i \subset S_C \setminus \Sigma_{S_C}$ with $D_i \subset D_{i+1}$ and $\bigcup_i D_i = S_C
\setminus \Sigma_{S_C}$. We consider the positive solutions of $L_\lambda \, f=0$ on the cone $C(D_i) \subset C$ over $D_i$ with vanishing boundary value along $\p C(D_i) \setminus \{0\}$.\\

Now we apply \cite[Remark 3.10]{L1} to  $L_\lambda$ on $C(D_i)$. Note that there exists $a_{C(D_i)}>0$ such that $dist(z,\Sigma_{S_C}) \le a_{C(D_i)} \cdot \delta_{\bp}(z) \le a_{C(D_i)} \cdot L_{\bp} \cdot  dist(z,\Sigma_{, aS_C})$. One readily checks the adaptedness conditions. Therefore, the Martin boundary of $L_\lambda$
equals $(\p C(D_i) \setminus \{0\}) \cup \{0_{C(D_i)},\infty_{C(D_i)}\}$ and the two fixed point solutions $\Psi_\pm[i]$ of $L_\lambda \,v$ vanishing along $\p C(D_i) \setminus \{0\}$, alternatively
labelled $0_{C(D_i)},\infty_{C(D_i)}$, are again of the form $\Psi_\pm[i](\omega,r) = \psi_\pm[i](\omega) \cdot r^{\alpha[i]_\pm}$.\\

Again, we insert $\Psi_\pm[i](\omega,r)$ into the equation $L_\lambda \, f =0$, written in polar coordinates as in \eqref{pol} and find that the $\psi_\pm[i]$ solve
\[
L_\lambda^\times \, v = - \Delta_{S_C} v(\omega) + \big(V^\times(\omega) - \lambda \cdot (\bp^\times)^2(\omega) \big)\cdot v(\omega) = (\alpha_\pm[i]^2 + (n-2)\cdot
\alpha_\pm[i]) \cdot v.
\]
Over $D_i$ we apply the spectral theory for bounded domains and observe that the \emph{positive} eigenfunctions $\psi_\pm[i]$ must equal the uniquely determined first Dirichlet eigenfunction $\psi[i]$ for the first eigenvalue $\mu[i]$ of $L_\lambda^\times$. Thus we have
\[
\psi[i]=\psi_-[i]=\psi_+[i] \,\mm{ and } \, \mu[i]=\alpha_-[i]^2 + (n-2)\cdot \alpha_-[i]=\alpha_+[i]^2 + (n-2)\cdot  \alpha_+[i],
\]
whence $\mu[i]\ge-(n-2)^2/4$. Further, the variational characterization of these eigenvalues gives $\mu[i]\ge \mu[i+1]$, since the space of admissible test functions on $D_i$ is a subset of the corresponding function space over $D_{i+1}$.\\

After a suitable normalization the first Dirichlet eigenfunctions $\psi[i]$ of $L_\lambda^\times$ on $D_i$ converge $C^3$-compactly on $S_C\setminus \Sigma_{S_C}$ as $i \ra \infty$ to a positive eigenfunction $\psi^*$ with eigenvalue $\mu^* = \lim_{i \ra \infty} \mu[i]\ge-(n-2)^2/4$. We also have the limits $\alpha_\pm^* := \lim_{i \ra \infty} \alpha_\pm[i]$. From this we observe as in \cite[Theorem 5.7]{L1} that $\psi^*>0$ is the non-weighted ground state of $L_\lambda^\times$ on $S_C \setminus \Sigma_{S_C}$ for the eigenvalue $\mu^* >-\infty$. Moreover, $L_\lambda^\times$, and hence $L^\times$, are shifted \si-adapted. Namely, $\bp^\times > c_{S_C} >0$ and thus the principal eigenvalue of $\delta^2_{\bp^\times} \cdot L_\lambda^\times$ remains finite as the non-weighted principal eigenvalue $\mu^*$ of $L_\lambda^\times$ is finite. Therefore, $L_\lambda^\times$ is a natural Schr\"odinger operator.\\

\textbf{Comparison Arguments} \,
Now we compare these solutions with $0_C,\infty_C \in \p_M (C,L_\lambda)$. From the fact that there are no positive solutions for $L_\lambda^\times \,  v = \mu \cdot v$ we deduce that if $\mu <\mu^*$,
\begin{equation}\label{aq}
\alpha_\pm^2 + (n-2)\cdot  \alpha_\pm\ge\mu^*\ge-(n-2)^2/4.
\end{equation}
In turn, the solution $\Psi_+=\psi(\omega) \cdot r^{\alpha_+}$ $L$-vanishes along $\widehat{\sigma_C} \setminus \{\infty\}$ and in particular in $0$. This implies the estimate
\begin{equation}\label{aqq}
\alpha_+^* \ge \alpha_+\mm{ and therefore  }\mu^* \ge  \alpha_\pm^2 + (n-2)\cdot  \alpha_\pm.
\end{equation}
From \eqref{aq} and \eqref{aqq} we conclude
\[
\mu^* =  \alpha_\pm^2 + (n-2)\cdot  \alpha_\pm \mm{ and therefore }
\alpha_\pm=\alpha^*_\pm.
\]
Since $0_C \neq \infty_C \in \p_M (C,L_\lambda)$, we further infer that
\[
\alpha_- <  -(n-2)/2 < \alpha_+ \mm{ and } \alpha_\pm^2 + (n-2)\cdot
\alpha_\pm > -(n-2)^2/4.
\]
Finally, the functions $\psi_\pm$ belong to the eigenvalue $\alpha_\pm^2 + (n-2)\cdot\alpha_\pm=\mu^*$. Since the ground state is uniquely determined, $\psi:=\psi_-=\psi_+=\psi^*$  is the ground state of $L_\lambda^\times$. \qed

%
\subsubsection{Geometric Operators on Cones} \label{ard}
Next we focus on a subclass of natural Schr\"odinger operators which one typically encounters in applications to scalar curvature geometry. For these we will derive quite explicit growth estimates.

\begin{theorem}[Eigenvalue Estimates for $J_C$ and $L_C$]\label{evee}
Let $C \in \mathcal{SC}_n$ be a singular area minimizing cone. Then we get the following estimates for the Jacobi field operator $J_C$ and the conformal Laplacian $L_C$:
\begin{itemize}
  \item The principal eigenvalue $\lambda_{J_C,C}$ of $J_C$ is non-negative and
  \[
  \mu_{C,(J_C)^\times_\lambda} \ge - \left(\frac{n-2}{2}\right)^2 \mm{ for any } \lambda \le \lambda_{J_C,C}.
  \]
  \item There are constants $\Lambda_n > \lambda_n >0$ depending only on $n$ such that $(L_C)_\lambda$ is \si-adapted for any $\lambda \le \Lambda_n$. Furthermore,
  \[
  \mu_{C,(L_C)^\times_\lambda} \ge - 1/4 \cdot \left(\frac{n-2}{2}\right)^2\mm{ for any }\lambda \le \lambda_n
  \]
  and thus
  \[\alpha_+ \ge - (1- \sqrt{3/4}) \cdot \frac{n-2}{2},\quad\alpha_- \le - (1+ \sqrt{3/4}) \cdot \frac{n-2}{2}\mm{ for }\lambda \le \lambda_n.
  \]
  Finally, we have $\mu_{C,(L_C)^\times_\lambda} <-\eta_\lambda <0$ for $\lambda \in (0,\lambda_n]$, whence
  \[
  -\vartheta_\lambda > \alpha_+ \ge - (1- \sqrt{3/4}) \cdot \frac{n-2}{2} >-\frac{n-2}{2} >- (1+ \sqrt{3/4}) \cdot \frac{n-2}{2} \ge \alpha_- > \vartheta_\lambda -(n-2)
  \]
  for constants $\eta_\lambda$ and $\vartheta_\lambda>0$ depending only on $\lambda$ and $n$.
\end{itemize}
\end{theorem}
\noindent\textbf{Proof} \,
From Proposition~\ref{scal} (iv) we know that $\lambda_{J_C,C} \ge 0$. For $\lambda <\lambda_{J_C,C}$, $J_C$ is \si-adapted and thus we have the lower bound $\mu_{C,(J_C)^\times_\lambda} >-(n-2)^2/4$. Moreover, for any given $\ve >0$, we find from Theorem~\ref{fix1} a common upper bound for $\mu_{C,(J_C)^\times_\lambda}$ for all $\lambda$ with $\lambda<  \lambda_{J_C,C} \le \lambda + \ve$.\\

Next we recall that the ground state $v>0$ is uniquely determined.  By standard elliptic theory applied to the given set of eigenvalue equations, $v$ is the unique $C^3$-limit of solutions $v_n >0$ of the shifted operator $(J_C)_{\lambda_{J_C,C}-1/n}\, f=0$, $n\geq1$. Note that in this case we have locally uniform bounds on the coefficients. Hence $\mu_{C,(J_C)^\times_{\lambda_{J_C,C}}} \ge - (n-2)^2/4$.\\

For $L_C$ we first notice that there is a common $\Lambda_n>0$ for all area minimizing cones in $\R^{n+1}$ such that $(L_C)_{\lambda}$ is \si-adapted for all $\lambda \le \Lambda_n$: For an individual cone this is just Proposition~\ref{scal} (i). The uniform estimate follows from the compactness of the space of singular cones. To proceed further, we note
\begin{itemize}
  \item the estimate  $\mu_{C,(J_C)^\times_0} \ge - \left(\frac{n-2}{2}\right)^2$.
  \item the Hardy inequality  for the \si-transform $\bp|_{S_C}$;  we can choose the constant $\tau(\bp,C)=\tau_n >0$ uniformly for all cones. This follows from the naturality of $\bp$.
  \item the identity $scal_C\equiv -|A|^2$ which is valid for minimal hypersurfaces in Ricci flat space, cf.\ the Gau\ss-Codazzi equation \eqref{miiq}.
\end{itemize}
This allows us to estimate the variational integral for $\mu_{C,(L_C)^\times_\lambda}$ for any smooth function $f$ compactly supported in $S_C \setminus \Sigma_{S_C} =\p B_1(0) \cap C \setminus \sigma_C$:
\begin{align*}
&\int_{S_C}  f  \cdot  (L_C)^\times_{\lambda} f  \,  dV \\
=& \int_{S_C} | \nabla f |^2 + \left(\frac{n-2}{4 (n-1)} \cdot scal_C|_{S_C} - \lambda \cdot \bp|_{S_C}^2\right)  \cdot  f^2 \, d V\\
\ge&\int_{S_C} | \nabla f |^2 - \frac{n-2}{4 (n-1)} \cdot |A|^2 \cdot f^2 dV - \max\{0,\lambda/\tau_n\} \cdot \int_{S_C}|\nabla f|^2  + |A|^2 \cdot f^2 dV.&
\end{align*}
For $\lambda/\tau_n < \frac{1}{4 (n-1)}$ and $|f|_{L^2}=1$ this yields the lower bound
\[
\ge \frac 14 \int_{S_C} | \nabla f |^2 -  |A|^2 \cdot f^2 dV \ge  \frac 14\cdot  \mu_{C,(J_C)^\times_0} \ge - \frac 14 \left(\frac{n-2}{2}\right)^2.
\]
Hence for $\lambda_n:=\min \{\Lambda_n, \frac{\tau_n}{4 (n-1)}\}/2$ we find
\[
\mu_{C,(L_C)^\times_\lambda} \ge - \frac 14\left(\frac{n-2}{2}\right)^2\mm{ for }\lambda \le \lambda_n.
\]
Next let $\lambda \in (0,\lambda_n]$. Since $scal_C \le 0$ and the codimension of $\Sigma_{S_C}$ is greater than two, we can find a function $f \in C_0^\infty(S_C \setminus \Sigma_{S_C})$ equal to $1$ outside a sufficiently small neighborhood of $\Sigma_{S_C}$, and equal to $0$ close to $\Sigma_{S_C}$ so that
\[
\int_{S_C} | \nabla f|^2 + \left(\frac{n-2}{4 (n-1)} \cdot scal_C|_{S_C} - \lambda \cdot \bp|_{S_C}^2\right)  \cdot  f^2 \, d V < -\frac 12 \int_{S_C} \lambda \cdot \bp|_{S_C}^2  \cdot  f^2 \, d V  < 0.
\]
Due to
\begin{itemize}
  \item the naturality of \si-transforms, and
  \item the \emph{compactness} of the space of all \emph{singular} area minimizing cones $C \subset \R^{n+1}$ in flat norm topology and in compact $C^5$-topology outside the singular sets,
\end{itemize}
we know that there is a common positive lower bound on $\int_{S_C}\bp|_{S_C}^2 \, d V$ for all such $C$. More precisely, there are constants $\eta_\lambda$, $\vartheta_\lambda>0$ which depend only on $\lambda>0$ and $n$ such that $\mu_{C,(L_C)^\times_\lambda} <-\eta_\lambda<0$, whence
\[
0>-\vartheta_\lambda>\alpha_+\quad\mm{and}\quad\alpha_- > \vartheta_\lambda -(n-2)>  -(n-2)
\]
for any singular area minimizing cone $C$. \qed\\

Inductive cone reduction arguments are a powerful tool to reduce general questions to low dimensions. In essence, this is an iterative blow-up process where we first blow-up around some $p_0 \in \Sigma_H$ and get a tangent cone $C_1$. While scaling around $0 \in \sigma_{C_1}$ merely reproduces $C_1$, blowing up around a singular point $p_1\neq 0 \in \sigma_{C_1}$ generates an area minimizing cone $C_2$ which can be written as a Riemannian product $C_2=\R \times C_3$ for a lower dimensional area minimizing cone $C_3$. This way we encounter cones $C^n\subset \R^{n+1}$ which can be written as a Riemannian product $\R^{n-k} \times C^k$, where $C^k \subset \R^{k+1}$ is a lower dimensional area minimizing cone singular at $0$. A famous instance of this technique is Federer's estimate for the codimension of the singularity set.\\

For applications to scalar curvature geometry which we shall consider elsewhere, it is useful to establish this kind of cone reduction arguments for the conformal Laplacian. In this case we observe that the minimal function $\Psi_+(\omega,r)$ of $(L_{C^n})_\lambda$ on $C^n$ shares the $\R^{n-k}$-translation symmetry with the underlying space $\R^{n-k} \times C^k$ since  $\Psi_+(\omega,r)$  is uniquely determined. Hence $\Psi_+(\omega,r)|_{\{0\} \times C^k}$ satisfies the equation
\[
(L_{C^k,n})_\lambda \Psi_+(\omega,r)|_{\{0\} \times C^k}=0,
\]
where
\[
L_{C^k,n}:=-\Delta - \frac{n-2}{4 (n-1)} \cdot |A|^2  \mm{ for }n \ge k \mm{ and } (L_{C^k,n})_\lambda  = L_{C^k,n}-\lambda  \cdot  \bp^2.
\]
Thus $L_{C^k,n}$ is a dimensionally shifted version of the conformal Laplacian on the cone $C^k$: The entities $\Delta$, $|A|^2$ and $\bp^2$ are intrinsically defined on $C^k$, while the dimensional shift comes from using $\frac{n-2}{4 (n-1)}$ in place of $\frac{k-2}{4 (k-1)}$. The next two results describe the analysis of these operators.

\begin{proposition}[Dimensionally Shifted $L_C$]\label{dimshift}
Let $C^k \subset \R^{k+1}$ be a singular area minimizing cone and $n \ge k$. Then we have:
\begin{itemize}
  \item $L_{C,n}$ is \si-adapted and for its principal eigenvalue there exists a uniform lower bound for all $k$-dimensional cones by a positive constant $\Lambda^*_k >0$ independent of $n$.
  \item There is a constant $\lambda^*_k \in (0,\Lambda^*_k)$ depending only on $k$ such that for any $n \ge k$ and $\lambda \le \lambda^*_k$,
  \[
  \mu_{C,(L_{C,n})^\times_\lambda} \ge - \frac{(k-2)^2}{12}
  \]
  and thus
  \[
  \alpha_+ \ge - (1- \sqrt{\frac 23}) \cdot \frac{k-2}{2},\quad \alpha_- \le - (1+ \sqrt{\frac 23}) \cdot \frac{k-2}{2}.
  \]
  Moreover, we have $\mu^*_{C,(L_{C,n})^\times_\lambda} <-\eta^*_\lambda <0$ for $\lambda \in (0,\lambda^*_k]$ and thus
  \[
  -\vartheta^*_\lambda >  \alpha_+ \ge - (1- \sqrt{\frac 23}) \cdot \frac{k-2}{2}>  - \frac{k-2}{2} >- (1+ \sqrt{\frac 23}) \cdot \frac{k-2}{2} \ge \alpha_- > \vartheta^*_\lambda -(k-2)
  \]
  for some constants $\eta^*_\lambda$ and $\vartheta^*_\lambda>0$ depending only on $\lambda$ and $k$.
\end{itemize}
\end{proposition}
\noindent\textbf{Proof} \,
From \eqref{kwsy} we get
\begin{align*}
\int_C | \nabla f |^2 - \frac{n-2}{4 (n-1)} \cdot |A|^2  \cdot  f^2 \, d V & \ge \int_C | \nabla f |^2 - \frac 14 \cdot |A|^2  \cdot  f^2 \, d V\\
&\ge\int_C \frac{k}{2 (k-1)} \cdot  |  \nabla f |^2 + \frac{k- 3}{4 (k-1)}\cdot   | A |^2   \cdot  f^2\,  d V,
\end{align*}
and for some $\tau^*_k >0$ depending only on the dimension this can be refined to
\[
\ge\frac{k-3}{4 (k-1)}\cdot \int_C f  \cdot  C_{H;A,0} f   \, d V \ge \tau^*_k  \cdot \frac{k- 3}{4 (k-1)}\cdot \int_C \bp^2 \cdot f^2   \, d V=:\Lambda^*_k \cdot \int_C \bp^2 \cdot f^2   \, d V
\]
where $C_{H;A,0}$ is the $A+B$-Laplacian from Chapter~\ref{exp}.\\

For the eigenvalue estimate on $S_C$ we write for any $\lambda < \Lambda^*_k$:
\begin{align*}
&\int_{S_C}  f  \cdot (L_{C,n})_\lambda f  \,  dV  = \int_{S_C} | \nabla f |^2 - \left(\frac{n-2}{4 (n-1)} \cdot | A |^2 +\lambda \cdot \bp^2\right)  \cdot  f^2 \, d V \ge\\ \ge&\int_{S_C} | \nabla f |^2 - \frac 14\cdot |A|^2 \cdot f^2 dV - \max\{0,\lambda/\tau_k\}  \cdot \int_{S_C}|\nabla f|^2  + |A|^2 \cdot f^2 dV.
\end{align*}
For $\lambda/\tau_k < 1/16$ and $|f|_{L^2}=1$, this can be refined to
\[
\ge \frac 13 \int_{S_C} | \nabla f |^2 -  |A|^2 \cdot f^2 dV \ge  \frac 13 \mu_{C,(J_C)^\times_0} \ge - \frac{(k-2)^2}{12}.
\]
Hence, for $\lambda^*_k:=\min \{\Lambda^*_k, 1/16\}/2$ we have
$\mu_{C,(L_{C,n})^\times_\lambda} \ge - \frac{(k-2)^2}{12}$ for $\lambda \le \lambda^*_k$.\\

For $\lambda \in (0,\lambda^*_k]$, we observe again, as in Proposition~\ref{evee}, that there is $f \in C_0^\infty(S_C \setminus \Sigma_{S_C})$ with
\[
\int_{S_C} | \nabla f|^2 + \left(\frac{n-2}{4 (n-1)} \cdot scal_C|_{S_C} - \lambda \cdot \bp|_{S_C}^2\right)  \cdot  f^2 \, d V < 0.
\]
Similarly as before, there exist constants $\eta^*_\lambda$ and $\vartheta^*_\lambda>0$, depending only on $\lambda>0$ and $k$, and such that $\mu^*_{C,(L_{C,n})^\times_\lambda} <-\eta^*_\lambda<0$. Therefore
\[
0>-\vartheta^*_\lambda>\alpha_+ \mm{ and }\alpha_- > \vartheta^*_\lambda -(k-2)>  -(k-2)
\]
for any singular area minimizing $k$-dimensional cone $C^k$. \qed\\

Complementary to the previous discussions where we mostly focussed on the radial growth rate, we now describe some global properties of the spherical component $\psi_C(\omega)$ of $\Psi_\pm[n,k](\omega,r) = \psi_C(\omega) \cdot r^{\alpha_\pm}$ which is defined over $S_C \setminus \Sigma_{S_C}$. In particular, this will yield uniform estimates for the radial growth near the singular set.

\begin{proposition}[Global Harnack Estimates for $\psi_C$]\label{skk}
For any $k$ with $n \ge k \ge 7$, $\lambda \in (0,\lambda^*_k]$ and cone $C^k \in {\cal{SC}}_k$, we consider the two solutions $\Psi_\pm[n,k](\omega,r) = \psi_C(\omega) \cdot r^{\alpha_\pm}$ of  $(L_{C^k,n})_\lambda \, f=0$ corresponding to $0$ resp.\ $\infty$ in the Martin boundary $\widehat{\sigma}_C$. Then there is a constant $a_{n,k,\lambda} >0$ depending on $n$, $k$ and $\lambda$, and a constant $b_{n,k,\lambda,\eta^*_\lambda}$ depending on $n$, $k$, $\lambda$ and $\rho>0$ defined by $\rho^2:=\lambda/\eta^*_\lambda$ (with $\eta^*_\lambda$ from Proposition~\ref{dimshift})
such that:
\begin{enumerate}
  \item $|\psi_C|_{L^p(S_C \setminus \Sigma_{S_C})}<\infty$ for $ p< \frac {k-1} {k-3}$.
  \item $|\psi_C|_{L^1(S_C \setminus \Sigma_{S_C})} \le a_{n,k,\lambda} \cdot \inf_{\omega \in S_C \setminus \Sigma_{S_C}}  \psi_C(\omega)$.
  \item $\sup_{\omega \in\E(\rho)}\psi_C(\omega) \le  b_{n,k,\lambda,\rho} \cdot|\psi_C|_{L^1(S_C \setminus \Sigma_{S_C})}$, where $\E(\rho)=\{x \in S_C \,|\, \bp(x) \le \rho^{-1}\}$.
\end{enumerate}
Similarly, if $q< \frac {k} {k-2}$, then $|v|_{L^q(B_1(0) \cap C^k)}<\infty$ for any solution $v>0$ of $(L_{C^k,n})_\lambda \, f=0$.
\end{proposition}
\noindent\textbf{Proof} \,
For $0 <\rho^2= \lambda/\eta^*_\lambda < -\lambda/\mu^*_{C,(L_{C,n})^\times_\lambda}$ and any $C \in {\cal{SH}}^{\R}_n$ we have
\begin{equation}\label{pp}
\frac{n-2}{4 (n-1)} \cdot |A|_C^2(\omega,1) + \lambda \cdot (\bp^\times)^2(\omega)  +\mu^*_{C,(L_{C,n})^\times_\lambda}>0 \mm{ on } \I(\rho).
\end{equation}
Hence $\psi_C(\omega)>0$ is a superharmonic function on $\I(\rho)$: $\psi_C$ solves the equation
\[
- \Delta_{S_C} \psi_C(\omega) - \Big(\frac{n-2}{4 (n-1)} \cdot |A|_C^2(\omega,1) + \lambda \cdot \bp^\times(\omega)^2 \Big)\cdot \psi_C(\omega) - \mu^*_{C,(L_{C,n})^\times_\lambda} \cdot \psi_C(\omega)=0.
\]
From this the pointwise estimate (\ref{pp}) implies $- \Delta_{S_C} \psi_C(\omega) \ge 0$  on $\I(\rho)$.\\

Although $S_C \subset \p B_1(0)$ is not a global area minimizer, it is an \emph{almost} minimizer and it shares the regularity theory with proper area minimizers. We can also locally apply the Bombieri-Giusti Harnack inequality \cite[Theorem 6 p.\ 39]{BG} in the following form: For any superharmonic function $w > 0$ defined on the regular region of $B_R(x) \cap S_C$ for a sufficiently small extrinsically measured radius $R>0$, we have
\[
0< \left\{\frac {1}{Vol_{n-1}(S_C \cap B_r(x))}\int w^p \right\}^{1/p} \leqq C  \cdot  \inf_{B_r(x)}w
\]
for $ r \leqq \beta_n \cdot R $ and $ p< \frac {k-1} {k-3}$ with constants $C = C(S_C,p)$ and $\beta_n > 0$. We apply this to a finite cover $B_r(p_j)$, $j=1,...m$, of $\Sigma_{S_C}$ by sufficiently small balls with $B_R(p_j) \subset \I(\rho)$.\\

Since the complement in $S_C \setminus \Sigma_{S_C}$ of these open balls is compact, $\inf_{S_C \setminus \Sigma_{S_C}}  \psi_C(\omega) >0$ and
$|\psi_C|_{L^p(S_C \setminus \Sigma_{S_C})} < \infty$. Moreover, we obtain for each individual cone $C$ the (trivial) estimate $\sup_{\omega \in\E(\rho)}\psi_C(\omega) \le  b\cdot|\psi_C|_{L^1(S_C \setminus \Sigma_{S_C})}$ for some suitably large $b=b_{C,n,k,\lambda,\rho}>0$. Since the $\psi_C$ are unique up to a multiple, the compactness of ${\cal{SC}}_k$, the naturality of $|A|$ and $\bp$, and the standard elliptic theory for $(L_{C,n})_\lambda^\times$ imply for all  $C^k \in {\cal{SC}}_k$ the existence of some common $a_{n,k,\lambda} >0$ such that
\[
|\psi_C|_{L^1(S_C \setminus \Sigma_{S_C})} \le a_{n,k,\lambda} \cdot \inf_{\omega \in S_C \setminus \Sigma_{S_C}}  \psi_C(\omega).
\]
In this way we also get a common $b_{n,k,\lambda,\rho}>0$ for all cones  such that  $\sup_{\omega \in\E(\rho)}\psi_C(\omega) \le b_{n,k,\lambda,\rho} \cdot|\psi_C|_{L^1(S_C \setminus \Sigma_{S_C})}$.\\

The assertion that $|v|_{L^q(B_1(0) \cap C^k)}<\infty$ for $q< \frac {k} {k-2}$ and a solution $v>0$ of $(L_{C^k,n})_\lambda \, f=0$ follows completely similarly by invoking again the Bombieri-Giusti Harnack inequality. \qed
%
%
%
\footnotesize
\renewcommand{\refname}{\fontsize{14}{0}\selectfont

\end{document}